\documentclass[11pt]{amsart}
\usepackage[utf8]{inputenc}
\usepackage[a4paper,width=185mm,top=20mm,bottom=15mm,bindingoffset=6mm]{geometry}
\usepackage[T1]{fontenc}
\usepackage{lmodern,amsmath,amssymb,amsthm,mathtools,mathrsfs}
\usepackage{enumitem,booktabs,array,microtype,needspace}
\usepackage{fancyhdr,tikz-cd}
\usepackage{float}
\usepackage{xcolor}

\usepackage[hidelinks]{hyperref}
\usepackage{aliascnt}
\newtheorem{theorem}{Theorem}[section]
\newaliascnt{lemma}{theorem}
\newtheorem{lemma}[lemma]{Lemma}
\aliascntresetthe{lemma}
\newaliascnt{proposition}{theorem}
\newtheorem{proposition}[proposition]{Proposition}
\aliascntresetthe{proposition}
\newaliascnt{corollary}{theorem}
\newtheorem{corollary}[corollary]{Corollary}
\aliascntresetthe{corollary}
\theoremstyle{definition}
\newtheorem{definition}{Definition}[section]
\newtheorem{remark}{Remark}[section]

\usepackage[nameinlink,noabbrev]{cleveref}
\crefname{example}{Example}{Examples}
\Crefname{example}{Example}{Examples}
\crefname{theorem}{Theorem}{Theorems}
\crefname{lemma}{Lemma}{Lemmas}
\crefname{proposition}{Proposition}{Propositions}
\crefname{corollary}{Corollary}{Corollaries}
\crefname{definition}{Definition}{Definitions}
\crefname{remark}{Remark}{Remarks}
\Crefname{theorem}{Theorem}{Theorems}
\Crefname{lemma}{Lemma}{Lemmas}
\Crefname{proposition}{Proposition}{Propositions}
\Crefname{definition}{Definition}{Definitions}
\Crefname{remark}{Remark}{Remarks}
\numberwithin{equation}{section}
\newcommand{\R}{\mathbb R}
\newcommand{\C}{\mathbb C}

\newcommand{\cR}{\mathcal{R}}
\newcommand{\cP}{\mathcal{P}}
\newcommand{\cB}{\mathcal{B}}

\newcommand{\dd}{\mathrm d}
\newcommand{\ii}{\mathrm i}
\setlist[enumerate]{leftmargin=*,itemsep=3pt,topsep=5pt}
\allowdisplaybreaks[2]
 
\makeatletter
\g@addto@macro\@adjustvertspacing{\advance\abovedisplayskip\z@ minus 2pt \advance\belowdisplayskip\z@ minus 2pt \advance\belowdisplayshortskip\z@ minus 2pt}
\makeatother
\normalsize
\title[Stratified superintegrability of spin CMS systems]{Stratified Superintegrability of Spin Calogero-Moser-Sutherland Systems}
\author[Guorui Ma]{Guorui Ma$^{1a,1b,*}$}
\author[Junze Zhang]{Junze Zhang$^{2}$}
\date{}
 
\makeatletter
\let\originalsetauthors\@setauthors
\renewcommand{\@setauthors}{%
  \originalsetauthors
  \begingroup
  \centering\footnotesize
  $^{1a}$ Shanghai Institute for Mathematics and Interdisciplinary Sciences (SIMIS), Shanghai 200433, China\par
  $^{1b}$ Research Institute of Intelligent Complex Systems, Fudan University, Shanghai 200433, China\par
  $^{2a}$ School of Mathematics and Physics, The University of Queensland, Brisbane, QLD 4072, Australia\par
   $^{2b}$ Department of Mathematical Sciences, Tsinghua University, Beijing, China \par
  \vspace{3pt}
  $^{*}$Corresponding author.\par
  \textit{Email:} \href{mailto:mgr18@tsinghua.org.cn}{\texttt{mgr18@tsinghua.org.cn}} (Guorui Ma)\quad
  \href{mailto:junze.zhang@uq.net.au}{\texttt{junze.zhang@uq.net.au}} (Junze Zhang)\par
  \endgroup
}
\makeatother
\hypersetup{pdftitle={Stratified Superintegrability of Spin Calogero–Moser–Sutherland Systems},pdfauthor={Guorui Ma and Junze Zhang}}
\begin{document}
\begin{abstract}
 We study the Poisson chain of the compact spin Calogero-Moser-Sutherland system obtained by conjugation reduction of $T^*\mathrm{SU}(n)$ at a fixed coadjoint orbit. We give a Horn-type condition ensuring that the spectral torus acts freely on the full fiber over a simple spectrum; for each orbit the admissible spectra form a nonempty open dense set, and the restricted chain is strongly superintegrable of rank $n-1$ on every orbit-type stratum. We also construct, across certain affine spectral walls, a proper quotient map onto a Hausdorff stratified Poisson space with a finite Whitney stratification. For $\mathrm{SU}(3)$, we compute the free spectral domain and the quotient near two intersecting walls, and a zero-orbit example shows that the analogous quotient need not be Hausdorff on the unrestricted reduced phase space. 
\end{abstract}
\maketitle

\section{Introduction}
\label{sec:introduction}

Conjugation reduction of $T^*\mathrm{SU}(n)$ at a fixed coadjoint
orbit gives a Poisson chain whose source may contain symplectic
strata of different dimensions \cite{KKS,Reshetikhin,MW,SL,ChenEtAl}.
We study two questions. Which feasible simple conserved spectra
have free spectral-torus action at every point of the corresponding
fiber, including its lower-dimensional orbit types? On the
pointwise free locus, when does subdivision by spectral walls admit
one Hausdorff quotient by the resulting characteristic leaves?

The first question leads to two spectral tests. Condition
\eqref{hyp:trace} excludes partial-trace equalities forced by a
nontrivial spectral stabilizer. Condition \eqref{hyp:spectrum}
tests whether the associated decompositions can be realized on both
complementary matrix blocks; it gives an exact criterion for freeness
on the complete fiber. For the second question, we construct the
characteristic quotient along selected integral affine walls and
verify its topology, smooth structure, and Poisson compatibility
across the refined strata \cite{LermanCuts,GJSImplosion,BL}.

\subsection{The model and the conserved spectrum}\label{subsec:main}
Fix $n\ge2$. Write $\ii=\sqrt{-1}$, let $A^*$ denote the conjugate transpose, and put $[n]=\{1,\ldots,n\}$. 
Fix $G=\mathrm{SU}(n)$, $\mathfrak{g}=\mathfrak{su}(n)$, and identify $\mathfrak{g}^*$ with $\mathfrak{g}$ by $\langle x,X\rangle=-\mathrm{tr}(xX)$. Choose
\begin{align*}
 \xi=\ii\mathrm{diag}(\mu_1,\ldots,\mu_n),\qquad  \mu_1\ge\cdots\ge\mu_n,\qquad \sum_i\mu_i=0,
\end{align*}
and let $\mathcal{O}_\xi=\{g\xi g^{-1}:g\in G\}$.

For a Hermitian matrix $A$, let $\mathrm{eig}(A)$ be its eigenvalues in nonincreasing order, counted with multiplicity, and let  $\lambda(x)=\mathrm{eig}(-\ii x)$. The ordered spectrum identifies $\mathfrak{g}/G$ homeomorphically with the closed Weyl chamber  \begin{align*}
 \mathcal{C}=\{\lambda\in\R^n:\lambda_1\ge\cdots\ge\lambda_n,\ \sum_i\lambda_i=0\}.
\end{align*}
The tuple $\lambda$ represents the orbit of $\ii\mathrm{diag}(\lambda_1,\ldots,\lambda_n)$. We write $\mathcal{C}^\circ$ for the subset with strict inequalities. A point $x$ with $\lambda(x)\in\mathcal{C}^\circ$, that is, with distinct eigenvalues, is said to have \emph{simple spectrum}. The smooth structures on the quotients below are described in Sections \ref{subsec:observable} and \ref{subsec:orbit-types}.

In the right trivialization $T^*G\cong\mathfrak{g}\times G$, the conjugation has a moment map. \begin{equation}
 J(x,\gamma)=x-\gamma^{-1}x\gamma. \label{eq:intro-model}
\end{equation}
Let $X_\xi =J^{-1}(\mathcal{O}_\xi)$. The reduced phase space is $\cR_\xi=X_\xi/G$. Here $/G$ denotes the conjugation orbit space, using simultaneous conjugation on pairs, and brackets denote an orbit. The two spaces of conserved quantities are \begin{align*} 
 \cP_\xi&=\{(x,y)\in\mathfrak{g} \times \mathfrak{g}:y\in\mathcal{O}_{-x},\ x+y\in\mathcal{O}_\xi\}/G,\qquad  \cB_\xi =\{\lambda(x):[(x,y)]\in\cP_\xi\}\subset\mathcal{C}.
\end{align*}
 The fixed-orbit Poisson chain of \cite[Eqs. (4.2)-(4.3)]{ChenEtAl} is
\begin{equation}\label{eq:poisschain}
 \cR_\xi\xrightarrow{\psi}\cP_\xi\xrightarrow{\pi}\cB_\xi,
 \qquad [(x,\gamma)]\longmapsto[(x,-\gamma^{-1}x\gamma)]\longmapsto\lambda(x).
\end{equation}
Write $y=-\gamma^{-1}x\gamma$, $z=J(x,\gamma)=x+y$, and $h=\pi\circ\psi$. Thus, $\cB_\xi$ is the set of ordered spectra compatible with the fixed orbit. We call a spectrum $\lambda\in\mathcal{C}$ \emph{feasible} if $\lambda\in\cB_\xi$, that is, if it is the spectrum of some point of $\cR_\xi$. Section \ref{subsec:observable} constructs the ambient chain before fixing $\mathcal{O}_\xi$ and explains this restriction.

We work with the connected symplectic orbit-type strata of the reduced space. Their definition, together with the associated quotient smooth structures, is recalled in Section \ref{sec:model}. Throughout, we distinguish the freeness of the spectral-torus action from the freeness of the conjugation action.

Let $\mathbf{H}=\{\mathrm{diag}(t_1,\ldots,t_n):|t_i|=1,\ \prod_i t_i=1\}$ be the diagonal maximal torus, and put $k=n-1$. We identify $\mathfrak h=\mathrm{Lie}(\mathbf{H})$ with the trace-zero real vectors by $b\mapsto\ii\mathrm{diag}b$, and identify $\mathfrak h^*$ with the same vector space using $\langle\lambda,b\rangle=\sum_i\lambda_i b_i$. Thus $h$ below is an $\mathfrak h^*$-valued moment map.

For a simple-spectrum $x$, choose $g\in G$ such that $g^{-1}xg=\ii\mathrm{diag}(\lambda(x))$.
For $q\in\mathbf{H}$, set $q(x)=gqg^{-1}$. This is independent of the ordered eigen-frame. The spectral torus acts by \begin{align*}
     q\cdot[(x,\gamma)]=[(x,q(x)\gamma)].
\end{align*} It is Hamiltonian on each symplectic orbit-type stratum, with
moment map $h$. The set of simple-spectrum points at which this
action is free is denoted by $\mathcal{M}_\xi^{\mathrm{free}}$. Put $\mathcal{P}_\xi^{\mathrm{free}}=\psi(\mathcal{M}_\xi^{\mathrm{free}})$
and $\mathcal{B}_\xi^{\mathrm{free}}=h(\mathcal{M}_\xi^{\mathrm{free}})$.

\subsection{Main results}
\label{subsec:main-result}
Condition \eqref{hyp:spectrum}, stated in Section \ref{subsec:spectralcond}, is a finite spectral test. It requires that, for every splitting of the index set, at least one
of two complementary Horn problems has no solution.
Define \begin{align}
\begin{split}
\Omega_\xi&=\{\lambda\in\mathcal{B}_\xi:\lambda\text{ satisfies }\eqref{hyp:spectrum}\},\\
\mathcal{M}_\xi&=h^{-1}(\Omega_\xi),\qquad
\mathcal{P}_\xi^\Omega=\pi^{-1}(\Omega_\xi).
\end{split}\label{eq:domains}
\end{align}
For any open restriction $U\subseteq\mathcal R_\xi$, its
\emph{retained strata} are the connected components of the
intersections of $U$ with the orbit-type strata.
Strong and stratified superintegrability are introduced in
Definitions \ref{def:strong} and \ref{def:stratified-system}.

\begin{theorem} 
\label{thm:main}
Let $\xi\in\mathfrak{su}(n)$ be arbitrary, define $\Omega_\xi$, $\mathcal{M}_\xi$ and $\mathcal{P}_\xi^\Omega$ by \eqref{eq:domains}, and equip the strata of $\mathcal R_\xi$ with their reduced symplectic forms \eqref{eq:reduction-form}. Then the following statements hold.
\begin{enumerate}[label=\textnormal{(\roman*)}]
\item For every feasible simple spectrum $\lambda$,
\begin{align*}
\lambda\in\Omega_\xi
&\Longleftrightarrow
G_{x,y}=G_{x,\gamma}\ \text{at every point of }h^{-1}(\lambda)
\\
&\Longleftrightarrow
\mathbf{H}\text{ acts freely on }h^{-1}(\lambda).
\end{align*}
The set $\Omega_\xi$ is a nonempty open semialgebraic subset of $\mathcal{C}^\circ$, dense in $\mathcal{B}_\xi$. Moreover, $\mathcal{M}_\xi$ is open and dense in $\mathcal{R}_\xi$. The spaces $\mathcal{M}_\xi$, $\mathcal{P}_\xi^\Omega$, and $\Omega_\xi$ have finite Whitney stratifications. On $\mathcal{M}_\xi$ use connected orbit-type strata, on $\mathcal{P}_\xi^\Omega$ use their images, and on $\Omega_\xi$ use connected components.
\item The sequence
\begin{align}\label{eq:main-diagram}
(\mathcal{M}_\xi,\omega_\xi)\xrightarrow{\ \psi\ }(\mathcal{P}_\xi^\Omega,\Pi_\xi)\xrightarrow{\ \pi\ }(\Omega_\xi,0)
\end{align}
 is stratified superintegrable in the sense of Definition \ref{def:stratified-system}.  
 All stratum closures are inherited by restriction to the open subset $\mathcal{M}_\xi\subset\mathcal R_\xi$: for $S\subset\mathcal{M}_\xi$, $\overline S^{\,\mathcal{M}_\xi}=\mathcal{M}_\xi\cap\overline S^{\,\mathcal R_\xi}$.
\item The map $\psi$ is the quotient by a free Hamiltonian action of $\mathbf{H}$. Locally it is a product with $\mathbf{H}$, including at singular strata. Both maps in \eqref{eq:main-diagram} are proper and surjective. If $\dim\mathcal{M}_\alpha=d_\alpha$ and $\mathcal{P}_\delta=\psi(\mathcal{M}_\alpha)$, put $\Pi_\delta=\Pi_\xi|_{\mathcal{P}_\delta}$. Then
\begin{align}\label{eq:main-dimensions}
\dim\mathcal{P}_\delta=d_\alpha-k,\qquad \mathrm{rank}\Pi_\delta=d_\alpha-2k,\qquad \mathrm{rank}d(\pi|_{\mathcal{P}_\delta})=k.
\end{align}
Each induced chain is strongly superintegrable. The connected components of the fibers of $\pi|_{\mathcal{P}_\delta}$ are its symplectic leaves.
\item The Hamiltonians $H_j(x)=\mathrm{tr}((-\ii x)^j)/j$, $2\le j\le n$, have complete flows. They preserve each stratum and $\mathcal{M}_\xi$. The same holds for a smooth function of these Hamiltonians on a neighborhood of their image on $\mathcal{M}_\xi$.
\end{enumerate}

\end{theorem}

The set $\Omega_\xi$ is maximal among simple spectra whose entire
fibers have spectral rank $n-1$ on every orbit-type stratum
(Remark \ref{rem:optimality}). The next theorem also applies to spectra $\lambda\notin\Omega_\xi$ for which the fiber $h^{-1}(\lambda)$ still contains points with trivial stabilizer for the spectral torus action.
Fixing linear functions $t_j(\lambda)=\langle\lambda,b_j\rangle-c_j$ selects spectral \emph{walls} $\{t_j=0\}$.
On the free locus over such a base, the theorem forms a quotient by identifying points related by the circle actions generated by those $t_j$ that vanish.

\begin{theorem}[Characteristic quotients along spectral walls]
\label{thm:integral-wall-quotient}
Let $b_1,\ldots,b_k$ be a basis of the integral lattice
$\Lambda=\{b\in\mathbb Z^n:\sum_i b_i=0\}$ of $\mathbf{H}$, with the convention
$b\mapsto\exp(\ii\theta\mathrm{diag}b)$ (period $2\pi$). Thus $k=n-1$.
Fix $1\le r\le k$, constants $c_1,\ldots,c_r\in\mathbb R$, and an open semialgebraic set
$B\subset\mathcal{B}_\xi^{\mathrm{free}}$. Set
\begin{align*}
t_j(\lambda)&=\langle\lambda,b_j\rangle-c_j,\qquad
M=\mathcal{M}_\xi^{\mathrm{free}}\cap h^{-1}(B),\qquad P=\psi(M).
\end{align*}
Decompose $B$ by the sign pattern of $(t_1,\ldots,t_r)$ and let $B_\beta$ be the connected components.

Let $M_\alpha$ be the connected components of the intersections of
$M$ with the original symplectic orbit-type strata, and put
$P_\delta=\psi(M_\alpha)$. Let $T$ and $S$ run through the connected
components of
\[
P_\delta\cap\pi^{-1}(B_\beta),\qquad
M_\alpha\cap\psi^{-1}(T),
\]
respectively. These are the refined middle and source strata.

For $\lambda\in B$, let the active-wall set be $I(\lambda)=\{j\in\{1,\ldots,r\}:t_j(\lambda)=0\}$ and put \begin{align*}
   K_I:=\exp\bigl(\mathrm{span}_{\mathbb R}\{\ii\mathrm{diag}b_j:j\in I\}\bigr)\subset\mathbf{H}. 
\end{align*} Define an equivalence relation on $M$ by
\begin{equation}\label{eq:wall-characteristic-relation}
m\sim m'\quad\Longleftrightarrow\quad m'\in K_{I(h(m))}\,m.
\end{equation}
Then the quotient map $p:M\to X:=M/{\sim}$ is proper and $X$ is Hausdorff.
Moreover, $X$ has a finite Whitney stratification (compatible with the above refinements) and admits compatible cone charts.

For every refined source stratum $S$, the map $p|_S:S\to p(S)$ is a
smooth characteristic reduction. There is a unique symplectic form
$\underline\omega_{p(S)}$ satisfying
\[
(p|_S)^*\underline\omega_{p(S)}=\omega|_S.
\]
The family $(\underline\omega_{p(S)})_S$ is Poisson compatible for
the smooth structure defined by the cone charts.

There is a proper smooth stratified Poisson map $\underline\psi:X\to P$ with $\psi=\underline\psi\circ p$. With these refinements and this characteristic quotient, the induced chain is stratified superintegrable in Definition \ref{def:stratified-system}.

\end{theorem}

We prove this theorem in Section \ref{sec:intewalquo}. 
The domains below satisfy
$\mathcal M_\xi\subseteq\mathcal M_\xi^{\mathrm{free}}$; the whole
domain of Theorem \ref{thm:main} lies in $M$ precisely when
$\Omega_\xi\subseteq B$.
\begin{center}
\small
\begin{tabular}{@{}p{0.26\linewidth}p{0.32\linewidth}p{0.32\linewidth}@{}}
\toprule
 & Theorem \ref{thm:main} & Theorem \ref{thm:integral-wall-quotient}\\
\midrule
Domain in $\mathcal R_\xi$ & $h^{-1}(\Omega_\xi)$ & $\mathcal M_\xi^{\mathrm{free}}\cap h^{-1}(B)$\\
Characteristic map & $p=\mathrm{id}$ & Torus collapse on the selected walls\\
Proper maps & $\psi$ and $\pi$ & $p$, $\psi$ and $\underline\psi$. $\pi$ need not be proper\\
Rank of the induced strong chain & $k$ & $k-c$ on a piece with $c$ active walls\\
\bottomrule
\end{tabular}
\end{center}
We use the Poisson-chain
approach of Chen, Jiang, Reshetikhin and Xiao. Remark \ref{rem:cjrx-conventions} states the
form conventions.

For $\xi=\ii\mathrm{diag}(a,0,-a)$ in $\mathfrak{su}(3)$,
Section \ref{sec:su3} computes both spectral domains and constructs
the quotient along $\lambda_1-\lambda_2=a$ and
$\lambda_2-\lambda_3=a$. Six- and four-dimensional orbit-type
strata occur in this construction. The zero-orbit example in
Proposition \ref{prop:unrestricted-zero} shows why the same leaf
relation cannot be used without a domain restriction.

\subsection{Relation with earlier work}
The theory of superintegrability goes back to Nekhoroshev
\cite{Nekhoroshev} and Mishchenko-Fomenko \cite{MF}.
The reductions of Kazhdan-Kostant-Sternberg \cite{KKS} and
Gibbons-Hermsen \cite{GH} produce Calogero systems and their
spin extensions. Reshetikhin \cite{Reshetikhin} describes the
conserved quantities by a Poisson chain. Global models for the
minimal orbit are given in \cite{FA,LMRS}.

Feh\'er \cite{Feher2024} proves constant-rank and integrability
results on suitable regular and principal domains. Feh\'er and
Fairon \cite[Lemmas 3.6 and 3.8 and Proposition 3.10]{FF2026}
construct the spectral torus and analyze its common stabilizers
with conjugation. Feh\'er \cite[Theorems 3.2 and 3.7]{Feher2026}
gives dimension formulas and generalized action-angle coordinates.
The discussion following \cite[Theorem 2.14]{FF2026} leaves
reductions outside the principal conjugation orbit type for
further study. Theorem \ref{thm:main} identifies complete free fibers
by complementary Horn problems and treats all orbit types
retained over these spectra, with their original incidences, that is,
with closure relations inherited from $\mathcal R_\xi$ in the restricted topology.

Singular Hamiltonian reduction provides the symplectic strata and
their local models \cite{SL,BL}. Chen, Jiang, Reshetikhin and Xiao
\cite{ChenEtAl} formulate stratified superintegrability using
Poisson chains and characteristic reduction.

Their Definition 3.5(7)-(8) requires the characteristic reductions to
form one stratified symplectic space, which is Hausdorff under their
Definition A.1, together with a quotient map through which the chain
factors. Separate smooth reductions do not by themselves establish
this global requirement. Sections
\ref{sec:unconditional}-\ref{app:disk} examine the obstruction caused
by nonclosed characteristic leaves in the unrestricted chain. Our positive results establish complete-fiber
freeness and construct a global characteristic quotient on the
stated restricted domains.

The local cone and the topological method in the wall theorem
have precedents. Clavijo Esquivel, Herbig and Seaton
\cite[Section 4]{HCES2024} study the same Poisson double cone.
Symplectic cuts collapse characteristic circles at a level set of the moment map
\cite{LermanCuts}. Guillemin, Jeffrey and Sjamaar
\cite[Lemma 2.3]{GJSImplosion} prove Hausdorffness of the imploded
quotient by showing that saturations of closed sets are closed.
Theorem \ref{thm:integral-wall-quotient} constructs the simultaneous
quotient along the selected walls in the fixed-orbit spin CMS system,
with compatible charts, Poisson brackets, and strong chains across
the orbit-type strata of its free domain.

\subsection{Organization of the paper}\label{subsec:organization}
Section \ref{sec:model} gives the reduction model and the definitions.
Section \ref{sec:proof} proves Theorem \ref{thm:main}. Section \ref{sec:intewalquo}
constructs the wall quotient. Section \ref{sec:su3} gives the
$\mathrm{SU}(3)$ calculation, the zero-orbit obstruction, and a second
obstruction on a compact $\mathrm{SU}(3)$ fiber.

\section{The reduction model and definitions}\label{sec:model}
 
We fix the symplectic signs, the quotient smooth structures, and
the meaning of the strong and stratified chains used in the two
main theorems.

\subsection{Cotangent geometry and conjugation reduction}\label{subsec:kks}
 
Write $I_r$ for the identity matrix, $\mathrm U(r)=\{V\in M_r(\C):V^*V=I_r\}$, and $S(\prod_j\mathrm U(r_j))$ for the block-diagonal subgroup with product of block determinants one. We identify $\mathfrak{g}\cong\mathfrak{g}^*$ using the invariant pairing
\begin{equation}\label{eq:pairing}
\langle x,X\rangle=-\mathrm{tr}(xX).
\end{equation}
In the right trivialization $T^*G\cong\mathfrak{g}\times G$, we fix the canonical one-form and sign conventions by
\begin{equation}\label{eq:canonical}
\Theta=-\mathrm{tr}(x\,\dd\gamma\gamma^{-1}),
\qquad \omega=\dd\Theta,
\qquad \iota_{X_f}\omega=-\dd f,
\qquad \{f,g\}=X_f(g).
\end{equation}
Represent a tangent vector by $(a,Z)$, with $\dot x=a$ and $\dot\gamma=Z\gamma$. Then
\begin{equation}\label{eq:tangent-form}
\omega((a,Z),(b,W))=-\mathrm{tr}(aW-bZ+x[Z,W]).
\end{equation}
Conjugation has moment map \eqref{eq:intro-model}, and the corresponding coadjoint-orbit form is
\begin{equation}\label{eq:orbit-form}
(\omega_{\mathcal{O}_\xi})_z([X,z],[Y,z])=\mathrm{tr}(z[X,Y]).
\end{equation}
The shifted Hamiltonian space is
\begin{equation}\label{eq:shift}
\mathscr M=T^*G\times\mathcal{O}_\xi^{-},\quad
\omega_{\mathscr M}=\omega-\omega_{\mathcal{O}_\xi},\quad
\Phi(x,\gamma,z)=J(x,\gamma)-z,
\end{equation}
where the superscript $-$ reverses the orbit form. Under the natural projection, $\Phi^{-1}(0)/G$ is identified with $\mathcal R_\xi$.

For a matrix $v$, set $G_v=\{g\in G:gv=vg\}$ and
\begin{equation}\label{eq:stabilizers}
G_{x,y}=G_x\cap G_y,\qquad G_{x,\gamma}=G_x\cap G_\gamma.
\end{equation}
For a closed subgroup $L\subset G$, the connected components of
$ (\mathcal R_\xi)_{(L)}=\{[(x,\gamma)]:G_{x,\gamma}\text{ is conjugate to }L\}$
are the original orbit-type strata. The center $Z_n=\{cI_n:c^n=1\}$ acts trivially, so the effective conjugation group is $G/Z_n$.
If $S$ is such a stratum, let $\widetilde S\subset X_\xi$ be its inverse image, $q:\widetilde S\to S$ the quotient, and $\iota:\widetilde S\hookrightarrow T^*G$ the inclusion. Its reduced symplectic form satisfies
\begin{equation}\label{eq:reduction-form}
q^*\omega_\xi=\iota^*\omega-(J|_{\widetilde S})^*\omega_{\mathcal{O}_\xi}.
\end{equation}
This is the usual reduction formula \cite{MW}, applied to \eqref{eq:shift}. For singular strata, use \cite[Theorems 7 and 15]{BL}. Compactness makes the action proper, without requiring $J$ to be proper. Condition \eqref{hyp:spectrum} will control the equality of the two stabilizers in \eqref{eq:stabilizers}, and hence the freeness of the spectral action, not of conjugation. 

\subsection{The ambient and fixed-orbit Poisson chains}\label{subsec:observable}
 
The map $(x,\gamma)\mapsto(x,-\gamma^{-1}x\gamma)$ has image
$\mathcal K=\{(x,y)\in\mathfrak{g} \times \mathfrak{g}:y\in\mathcal{O}_{-x}\}$ and induces the ambient chain \cite[Eqs. (4.5)-(4.7)]{ChenEtAl}
\begin{equation}\label{eq:intro-ambient-chain}
T^*G/G\xrightarrow{\widetilde\psi}\mathcal K/G
\xrightarrow{\widetilde\pi}\mathfrak{g}/G\cong\mathcal{C},
\qquad [(x,\gamma)]\mapsto[(x,-\gamma^{-1}x\gamma)]\mapsto\lambda(x).
\end{equation}
All quotients below are taken with respect to the diagonal $G$-action by conjugation.
We use the usual smooth structure on these quotients: a function is smooth if, locally, it is induced by a $G$-invariant smooth function on the ambient manifold. In particular, a function $f$ on $\mathcal{C}\cong\mathfrak{g}/G$ is smooth exactly when $f\circ\lambda$ is smooth on $\mathfrak{g}$.
Imposing
\begin{equation}\label{eq:original-spaces}
\widetilde{\mathcal{P}}_\xi=\{(x,y)\in\mathcal K:x+y\in\mathcal{O}_\xi\}
\end{equation}
gives the fixed-orbit chain \eqref{eq:poisschain}. Both maps are surjective: every pair in \eqref{eq:original-spaces} has a lift $y=-\gamma^{-1}x\gamma$ in $G$, whose moment map is $x+y$.

The two momentum factors commute, and their product bracket, with gradients for \eqref{eq:pairing}, is
\begin{equation}\label{eq:product-poisson}
\{f,g\}_{\rm prod}=\mathrm{tr}(x[\nabla_xf,\nabla_xg])+\mathrm{tr}(y[\nabla_yf,\nabla_yg]).
\end{equation}
Its Hamiltonian flows preserve the spectra of $x,y$. For invariant Hamiltonians they also preserve the diagonal moment map $x+y$. Thus, the ideal of invariant functions vanishing on \eqref{eq:original-spaces} is a Poisson ideal, and the bracket descends independently of invariant local extensions. Reduction gives
\begin{equation}\label{eq:original-poisson}
\{f\circ\psi,g\circ\psi\}_{\omega_\xi}=\{f,g\}_{\mathcal{P}_\xi}\circ\psi.
\end{equation}
Invariant functions of $x$ are Casimirs on the middle space, so $\pi$ is Poisson for the zero bracket on $\mathcal{B}_\xi$. Lemma \ref{lem:global-torus-quotient} identifies its restricted smooth structure with the spectral-torus quotient structure.

The Horn theorem \cite[Theorem 1]{Fulton}, applied to spectra $\lambda$, $(-\lambda_n,\ldots,-\lambda_1)$ and $\mu$ of $-\ii x,-\ii y,-\ii(x+y)$, describes $\mathcal{B}_\xi$ by finitely many affine inequalities in $\lambda$. Conversely, a Hermitian realization lifts to $G$ by adjusting a unitary conjugating matrix by a scalar. Hence, $\mathcal{B}_\xi$ is a convex polyhedron. Its faces need not be images of orbit-type strata, and the components of $\Omega_\xi$ need not be convex.

\subsection{Particle coordinates and Hamiltonian functions}\label{subsec:particle}

Let $G_{\rm reg}$ be the elements of $G$ with distinct eigenvalues and $\mathbf{H}_{\rm reg}=G_{\rm reg}\cap\mathbf{H}$. Diagonalizing $\gamma$ on this locus gives the standard particle model $(T^*\mathbf{H}_{\rm reg}\times\mathcal{O}_\xi')/\mathbf{W}$, where $\mathcal{O}_\xi'=\{z\in\mathcal{O}_\xi:z_{ii}=0\}/\mathbf{H}$ and $\mathbf{W}\cong\mathfrak S_n$ acts by simultaneous permutations \cite{Reshetikhin,ChenEtAl}. The main theorems impose no regularity
on $\gamma$; we use the invariant Hamiltonians and flows below.

For $2\le j\le n$, define the real Hamiltonians and their skew-Hermitian gradients by
\begin{align}
 H_j(x)&=\frac1j\mathrm{tr}((-\ii x)^j), \label{eq:gradients}\\
 D_j(x)&=\ii\left((-\ii x)^{j-1}-\frac{\mathrm{tr}((-\ii x)^{j-1})}{n}I_n\right),\qquad \dd H_j(a)=-\mathrm{tr}(D_j(x)a).\nonumber
\end{align}
Since $[D_j(x),x]=0$, \eqref{eq:tangent-form} gives the complete Hamiltonian flow
\begin{align}\label{eq:original-flow}
 x(t)=x(0),\qquad \gamma(t)=\exp(tD_j(x(0)))\gamma(0).
\end{align}
The exponential lies in $G$. The matrices $y,z$ and the functions $H_l$ remain fixed. Hence, the pullback of every smooth function on $\cP_\xi$ is a first integral.

\subsection{Strong Poisson chains}
 The reductions of each stratum will satisfy the following definition, as stated in \cite[Definition 3.1]{ChenEtAl}. 

\begin{definition}[Strong superintegrability]\label{def:strong}
 Let $(M,\omega)$ be a symplectic manifold of dimension $2r$, and let $0\le k\le r$. A strong superintegrable system of rank $k$ consists of Poisson maps \begin{align}
 (M,\omega)\xrightarrow{\psi}(P,\Pi_P)\xrightarrow{\pi}(B,0). \label{eq:strong}
\end{align}
Here $(P,\Pi_P)$ is a Poisson manifold, and $B$ is a smooth manifold with zero Poisson bracket, with $\dim P=2r-k$ and $\dim B\ge k$. The map $\psi$ is a surjective submersion, and $\pi$ has constant rank $k$. The connected components of every nonempty fiber of $\pi$ are the symplectic leaves of $P$.
\end{definition}

 The case $k=r$ is Liouville integrability, and when $0<k<r$ the system is strictly superintegrable. Rank zero is allowed. Completeness is a separate requirement, and $\pi$ need not be surjective. For $f\in C^\infty(B)$ and $g\in C^\infty(P)$, the Casimir identity gives $\{(\pi\circ\psi)^*f,\psi^*g\}=0$.

\subsection{Stratified spaces} \label{sec:stratification}

We use semialgebraic sets in the usual sense of finite
unions of sets defined by polynomial equalities and inequalities
\cite[Section 2.1]{BCR}.

We use local ambient smooth extensions on singular spaces, with the following conventions.

\begin{definition}\label{def:whitney-space}
A stratification of a Hausdorff paracompact space $X$ is a locally finite decomposition into connected, locally closed smooth manifolds without boundary, satisfying the frontier condition
\begin{align}\label{eq:frontier}
 X_\alpha\cap\overline{X_\beta}\ne\varnothing
 \quad\Longrightarrow\quad X_\alpha\subseteq\overline{X_\beta}.
\end{align}
A singular smooth atlas identifies open subsets of $X$ with locally closed Euclidean subsets, with smooth embeddings on strata and local ambient smooth extensions of all transitions and their inverses. The sheaf $C_X^\infty$ consists of functions with local ambient smooth extensions. A continuous map $F:X\to Y$ is smooth if $F^*C_Y^\infty(U)\subseteq C_X^\infty(F^{-1}(U))$ for every open $U\subseteq Y$, and stratified if it maps each stratum into one stratum. An isomorphism is a bijection whose map and inverse are smooth and stratified.

The stratification is Whitney if, in these charts, for $x_j\in X_\alpha$ tending to $y\in X_\beta\subseteq\overline{X_\alpha}$ with $T_{x_j}X_\alpha\to T$, one has $T_yX_\beta\subseteq T$. If $y_j\in X_\beta$ also tends to $y$ with $x_j\ne y_j$ and $\mathbb R(x_j-y_j)\to\ell$, then $\ell\subseteq T$. These are conditions (A) and (B). It is finite if it has finitely many strata.
\end{definition}

\subsection{Hamiltonian quotients and reduced two-forms}\label{subsec:orbit-types}

For a Hamiltonian action of a compact Lie group $K$ on a symplectic manifold $(N,\sigma)$ with equivariant moment map $\Phi:N\to\mathfrak k^*$, where $\mathfrak k=\mathrm{Lie}(K)$, put $Z=\Phi^{-1}(0)$ and $q:Z\to X=Z/K$. By singular reduction, $X$ decomposes into symplectic orbit-type strata, and each stratum carries the reduced form appearing in \eqref{eq:reduction-form} \cite[Theorems 7 and 15]{BL}. If $f\in C^\infty(X)$, then $q^*f$ is $K$-invariant on $Z$ and extends locally to a $K$-invariant smooth function $\widetilde f$ on $N$ (and similarly for $g$). For any such choices $\widetilde f,\widetilde g$, the reduced Poisson bracket on $X$ is characterized by
\begin{equation}\label{eq:reduced-bracket}
q^*\{f,g\}_X=\{\widetilde f,\widetilde g\}_N|_Z.
\end{equation}
It is independent of the extensions and restricts to the symplectic bracket on each stratum. On submanifolds of strata, we use the restricted reduced forms.

\begin{definition} \label{def:stratified-poisson}
Let $X$ be a Whitney stratified space equipped with a singular smooth atlas.
A \emph{stratified Poisson structure} on $X$ is a Poisson bracket on the smooth functions $C^\infty(X)$ such that, for each stratum $S\subset X$, the restricted bracket on $C^\infty(S)$ is the usual Poisson bracket of a smooth bivector field $\Pi_S\in\Gamma(\bigwedge^2 TS)$.
Moreover, this bivector is locally induced from an ambient smooth bivector: for every point of $X$ and every compatible singular chart $\iota:U\hookrightarrow\mathbb R^r$, there exists a smooth bivector field $\widetilde\Pi$ on an open neighborhood of $\iota(U)$ in $\mathbb R^r$ whose restriction to each stratum agrees with $\Pi_S$ under the chart.  See \cite[Definition 2.1]{ChenEtAl}.

A family of symplectic forms on the strata of $X$ is \emph{Poisson compatible} if their induced Poisson brackets arise as the restrictions of some stratified Poisson bracket on $C^\infty(X)$.
\end{definition}

\begin{definition} \label{def:refine}
A Whitney stratification $\{S_\gamma\}$ of $X$ refines $\{X_\alpha\}$ if each $S_\gamma$ is a smooth embedded submanifold of some $X_\alpha$.
For a stratified Poisson space $(X,\Pi)$, the refinement is \emph{Poisson} if $\Pi(x)\in\bigwedge^2T_xS_\gamma$ for every $x\in S_\gamma$. For an open subset $U\subseteq X$, the restricted stratification consists of the connected components of $U\cap X_\alpha$.
\end{definition}

\subsection{Stratified superintegrability and characteristic quotients}
\label{subsec:stratified-definition}
We now introduce the notion of a \emph{stratified superintegrable chain}. The main issue is that the relevant two-forms are typically only presymplectic on refined strata, so we record the resulting characteristic leaf-space quotients and how they fit together globally.

\begin{definition} 
\label{def:characteristic-reduction}
A \emph{presymplectic form} on a manifold $Y$ is a closed two-form
$\tau$ of constant rank. Its kernel distribution is
$\ker\tau_y=\{v\in T_yY:\tau_y(v,w)=0\text{ for all }w\in T_yY\}$.
Its maximal connected integral manifolds are the \emph{characteristic
leaves}. Write $y\sim_\tau y'$ when both points lie on the same leaf.
A \emph{smooth characteristic reduction} is the leaf space
$Q=Y/{\sim_\tau}$ with its quotient topology and a smooth manifold
structure for which the projection $p:Y\to Q$ is a submersion.
\end{definition}
Since $\tau$ is closed, the distribution $\ker\tau$ is integrable. For a smooth characteristic reduction $p:Y\to Q$, one has $\ker dp=\ker\tau$ and connected fibers, and $\tau$ is the pullback of a unique symplectic form $\omega_Q$ on $Q$ with $p^*\omega_Q=\tau$.

\begin{definition}[Stratified superintegrability]\label{def:stratified-system}
A stratified superintegrable chain consists of the following data.
All stratified spaces below are Hausdorff and have finite Whitney
stratifications and singular smooth atlases as in
Definition \ref{def:whitney-space}.

\emph{Original chain.} The maps
$\mathcal{M}\xrightarrow{\psi}\mathcal{P}\xrightarrow{\pi}\mathcal{B}$
are smooth and Poisson, and $\psi$ is surjective. The strata
$\mathcal{M}_\alpha$ have Poisson-compatible symplectic forms
$\omega_\alpha$. The space $(\mathcal{P},\Pi_{\mathcal{P}})$ is
stratified Poisson, with strata $\mathcal{P}_\delta$, and
$\mathcal{B}$ has zero bracket and strata $\mathcal{B}_\beta$.

\emph{Compatible refinements.} The connected components $T$ of
$\mathcal{P}_\delta\cap\pi^{-1}(\mathcal{B}_\beta)$ and the connected
components $S$ of $\mathcal{M}_\alpha\cap\psi^{-1}(T)$ form finite
Whitney refinements. The refinement $\{T\}$ is Poisson, with
$\Pi_T=\Pi_{\mathcal{P}}|_T$, and
$\tau_S=\omega_\alpha|_S$ is presymplectic. We do not assume that the Poisson bivector on $\mathcal{M}_\alpha$ is tangent to $S$.

\emph{Global characteristic quotient.} Define
\begin{equation}\label{eq:separate-relation}
m\sim m'\quad\Longleftrightarrow\quad
m,m'\in S\text{ for some }S\text{ and }m\sim_{\tau_S}m'.
\end{equation}
The projection $p:\mathcal{M}\to\underline{\mathcal{M}}=\mathcal{M}/{\sim}$
has the quotient topology on its target and is smooth and stratified
for the refinement $\{S\}$. The strata of $\underline{\mathcal{M}}$
are $Q_S=p(S)$, and $p_S=p|_S:S\to Q_S$ is a smooth characteristic
reduction. Its symplectic forms are Poisson compatible and satisfy
\begin{equation}\label{eq:global-characteristic-pullback}
p_S^*\underline\omega_{Q_S}=\tau_S.
\end{equation}
Write $\omega^{\rm ref}=(\tau_S)_S$ and
$\underline\omega=(\underline\omega_{Q_S})_S$.
A smooth stratified Poisson map $\underline\psi$ makes the diagram
\begin{equation}\label{eq:global-reduction-diagram}
\begin{tikzcd}[column sep=large,row sep=large]
\mathcal{M} \arrow[r,"\psi"] \arrow[d,"p"'] &
\mathcal{P} \arrow[r,"\pi"] \arrow[d,equal] &
\mathcal{B} \arrow[d,equal] \\
\underline{\mathcal{M}} \arrow[r,"\underline\psi"'] &
\mathcal{P} \arrow[r,"\pi"'] & \mathcal{B}
\end{tikzcd}
\end{equation}
commute, where both copies of $\mathcal{P}$ carry the refinement $\{T\}$.
Thus $p$ collapses precisely the characteristic leaves within each $S$.

\emph{Strong chains.} For each $S$ over $T$ and $\mathcal{B}_\beta$,
the bottom row of the commutative diagram
\begin{equation}\label{eq:stratum-reduction-diagram}
\begin{tikzcd}[column sep=large,row sep=large]
(S,\tau_S) \arrow[r,"\psi|_S"] \arrow[d,"p_S"'] &
(T,\Pi_T) \arrow[r,"\pi|_T"] \arrow[d,equal] &
(\mathcal{B}_\beta,0) \arrow[d,equal] \\
(Q_S,\underline\omega_{Q_S}) \arrow[r,"\underline\psi|_{Q_S}"'] &
(T,\Pi_T) \arrow[r,"\pi|_T"'] & (\mathcal{B}_\beta,0)
\end{tikzcd}
\end{equation}
is strong in the sense of Definition \ref{def:strong}.
The top-left space is presymplectic, and the bottom-left space is its
symplectic reduction. Only the bottom row is required to be a strong
Poisson chain.
\end{definition}

For $\underline{\mathcal{M}}$ to be Hausdorff, each characteristic leaf must be closed in $\mathcal{M}$, not just within its stratum $S$.

\begin{remark} \label{rem:cjrx-conventions}
Definition \ref{def:stratified-system} requires the reduced forms to induce a stratified Poisson bracket on $C^\infty_{\underline{\mathcal{M}}}$. It does not require them to extend to an ambient two-form. In Theorem \ref{thm:main}, $p=\mathrm{id}$ and $\psi$ is the spectral torus quotient. On the walls of Theorem \ref{thm:integral-wall-quotient}, $p$ collapses characteristic subtori and $\underline\psi$ is the
quotient by the residual torus.

\end{remark}

\section{Complete spectral fibers}\label{sec:proof}
 
The proof of Theorem \ref{thm:main} proceeds in three parts. We begin by relating the complementary Horn condition to the freeness of the spectral torus action: the condition fails precisely when the torus has a nontrivial stabilizer at some point in the fiber. Next, we show that the spectra satisfying the condition form a dense subset of the feasible region. In the last step, we identify the torus quotient with the middle space $\mathcal{P}_\xi^\Omega$ and verify the stated rank and leaf properties on each retained stratum. 

\subsection{The complementary Horn condition}
\label{subsec:spectralcond}
\ Fix an ordered spectrum $\lambda = (\lambda_1,\ldots,\lambda_n)$ and a fixed orbit spectrum $\mu=(\mu_1,\ldots,\mu_n)$, both written in non-increasing order. For any subset $I=\{i_1<\cdots<i_r\}\subset[n]$, define its complement $I^c=[n]\setminus I$, the corresponding sub-tuple $\lambda[I]:=(\lambda_{i_1},\ldots,\lambda_{i_r})$, and the partial sum $\lambda_I:=\sum_{i\in I}\lambda_i$. We use the analogous notation $\mu[I]$ and $\mu_I$ for $\mu$. We regard subsets as indexing \emph{positions} in the ordered lists (so repetitions in $\mu$ are counted with multiplicity), and we do not assume that different subsets are disjoint. Finally, for any real tuple $a$, we write $a^\downarrow$ for its non-increasing rearrangement. In particular, $(-\lambda[I])^\downarrow=(-\lambda_{i_r},\ldots,-\lambda_{i_1})$.

We seek a condition on a feasible simple spectrum that excludes nontrivial spectral-torus stabilizers at every point of its fiber. Proposition \ref{prop:stabilizer} detects such a stabilizer by a nonzero proper common reducing subspace $W$ of $x$ and $y$ with $\gamma W\ne W$. The simple spectrum of $x$ gives distinct index sets $I,J$ of the same cardinality, corresponding to $W$ and $\gamma W$. If $-\ii(x+y)|_W$ has spectrum $\mu[K]$, taking traces gives \begin{align*}
    \lambda_I-\lambda_J-\mu_K=0.
\end{align*} 
The restrictions on $W$ and $W^\perp$ must also solve both complementary matrix-sum problems. We can therefore exclude every such trace equality or test whether the two block problems are simultaneously realizable. The first option gives the simpler
sufficient condition \begin{align}
&\lambda\in\mathcal{C}^\circ,\quad \lambda_I-\lambda_J\ne\mu_K, \   \text{for every integer }r \in\{1,\ldots,n-1\} \notag\\
&\text{ and }I,J,K\subset[n]\text{ with }|I|=|J|=|K|=r,\ I\ne J.\tag{\ensuremath{H^{\prime}}}\label{hyp:trace}
\end{align}
Proposition \ref{prop:comparison} proves that \eqref{hyp:trace} implies the condition \eqref{hyp:spectrum} below and provides a spectrum that satisfies \eqref{hyp:spectrum} but not \eqref{hyp:trace}.

 The weaker condition is stated in terms of the following set of pairs of Hermitian matrices. 
\begin{definition}
\label{def:Horn}
Let $a,b,c$ be real $m$-tuples in non-increasing order, where $m\ge1$, and let $\mathrm{Herm}_m(\C)$ be the space of Hermitian $m\times m$ matrices. Define
\begin{align*}
\mathrm{Horn}_m(a,b;c) &:= \bigl\{(A,B)\in \mathrm{Herm}_m(\C)\times\mathrm{Herm}_m(\C):\ \mathrm{eig}(A)=a,\ \mathrm{eig}(B)=b,\ \mathrm{eig}(A+B)=c\bigr\}.
\end{align*}
We say that the \emph{Horn problem} for $(a,b;c)$ is solvable when this set is nonempty. We call $\mathrm{Horn}_m(a,b;c)$ the \emph{Horn set}. For $m=1$, this reduces to the scalar identity $c_1=a_1+b_1$.
\end{definition}

 By the Horn theorem \cite[Theorem 1]{Fulton}, $\mathrm{Horn}_m(a,b;c)\ne\varnothing$ if and only if $\sum_i c_i=\sum_i a_i+\sum_i b_i$ and finitely many linear inequalities in $a,b,c$ hold. We use this finite description and the matrix formulation above.  

The spectral hypothesis used in Theorem \ref{thm:main} is
\begin{align}
\tag{H}\label{hyp:spectrum}
&\lambda\in\mathcal{C}^\circ, \ \text{for every integer }r\in\{1,\ldots,n-1\}\notag\\
& \text{ and every }I,J,K\subset[n]\text{ with }|I|=|J|=|K|=r\text{ and }I\ne J,\notag\\
&\mathrm{Horn}_r\bigl(\lambda[I],(-\lambda[J])^\downarrow;\mu[K]\bigr)=\varnothing \ \text{or}\  \mathrm{Horn}_{n-r}\bigl(\lambda[I^c],(-\lambda[J^c])^\downarrow;\mu[K^c]\bigr)=\varnothing. \notag
\end{align} Here $\varnothing$ means that an obstructing decomposition is infeasible. For a feasible simple spectrum, Propositions \ref{prop:stabilizer} and \ref{prop:comparison} give \begin{align*}
    \eqref{hyp:trace}\ \Longrightarrow\ \eqref{hyp:spectrum} \ \Longleftrightarrow\ \mathbf{H}\text{ acts freely at every point of }h^{-1}(\lambda).
\end{align*} These are alternative tests, rather than two assumptions to impose together. The trace test also supplies the density argument, while
the complementary Horn test characterizes the complete free fibers.

The two Horn problems correspond to independent matrix sums on complementary nonzero blocks. Taking the trace in the $r\times r$ problem forces the identity $\lambda_I-\lambda_J=\mu_K$. Thus, the stronger trace condition \eqref{hyp:trace} implies \eqref{hyp:spectrum}. Condition \eqref{hyp:spectrum} is weaker in the sense that it permits $\lambda_I-\lambda_J = \mu_K$, provided that at least one of the two block Horn problems has no solution.

For $S\subset\mathcal{B}_\xi$, restricting spectra in $S$ just means: keep all points in the fibers lying over $S$, and use the same structure maps and the subspace topology inherited from the original spaces.

\subsection{The spectral condition and the stabilizers}

In the rest of this paper, take $\mathfrak{g}_{\mathrm{reg}}=\{x\in\mathfrak{g}:\lambda(x)\in\mathcal{C}^\circ\}$. For such $x$, and for each $i \in [n]$, let $E_i(x)$ be the orthogonal projection onto its $\ii\lambda_i(x)$-eigenspace. For $I\subset[n]$, write $E_I(x)=\bigoplus_{i\in I}\mathrm{Im}E_i(x)$. For $m=(x,\gamma)$, put $y=-\gamma^{-1}x\gamma$. The inclusion $G_{x,\gamma}\subset G_{x,y}$ always holds.  The next proposition states when equality holds at every point of the whole fiber $h^{-1}(\lambda)$. 

\begin{proposition}\label{prop:stabilizer}
Let $\lambda\in\mathcal{B}_\xi$ be feasible (Section \ref{subsec:main}) with distinct entries. The following conditions are equivalent:
\begin{enumerate}
\item The spectrum $\lambda$ satisfies \eqref{hyp:spectrum}.
\item For every $(x,\gamma)$ with $\mathrm{eig}(-\ii x)=\lambda$ and $J(x,\gamma)\in\mathcal{O}_\xi$, we have $G_{x,y}=G_{x,\gamma}$.
\item The action of $\mathbf{H}$ defined in \eqref{eq:eigenvalue-action} below is free on $h^{-1}(\lambda)$.
\end{enumerate}
\end{proposition}

\begin{proof}
We first relate (2) and (3) by computing the stabilizer of the $\mathbf{H}$-action on a fiber $h^{-1}(\lambda)$. Then we show that (1) holds if and only if the equality $G_{x,y}=G_{x,\gamma}$ in (2) holds for every point in that fiber.

Assume that all the eigenvalues $\lambda_1,\ldots,\lambda_n$ of $-\ii x$ are distinct (ordered as in Section \ref{subsec:main}), and set \begin{align*}
E_i(x)=\prod_{j\neq i}\frac{-\ii x-\lambda_j I_n}{\lambda_i-\lambda_j},\qquad i=1,\ldots,n.
\end{align*}
Then $xE_i(x)=\ii\lambda_iE_i(x)$ and $\sum_{i=1}^n E_i(x)=I_n$. In terms of these matrices, the spectral action of Section \ref{subsec:main} reads, for $t=\mathrm{diag}(t_1,\ldots,t_n)\in\mathbf{H}$,
\begin{align}\label{eq:eigenvalue-action} t(x)=\sum_{i=1}^n t_iE_i(x),\qquad t\cdot(x,\gamma)=(x,t(x)\gamma). \end{align}
On $\mathfrak{g}_{\mathrm{reg}}$, the eigenprojections $E_i(x)$ vary smoothly with $x$. The $\mathbf{H}$-action also commutes with simultaneous conjugation and leaves $x$, $y$, and $J$ unchanged. Therefore, it restricts to an action on $h^{-1}(\mathcal{C}^\circ)\subseteq\mathcal R_\xi$. Let $T_x:=G_x$ (an abelian group). If two points in $\mathcal R_\xi$ project to the same pair $(x,y)$, then their $\gamma$'s differ by left-multiplication by some $t\in T_x$, i.e. the lifts are of the form $(x,t\gamma)$. Now fix $g\in G_{x,y}$. Then both $g$ and $\gamma g\gamma^{-1}$ lie in $T_x$, so \begin{align}
\chi_\gamma:G_{x,y}\longrightarrow T_x,\qquad \chi_\gamma(g) = g\gamma g^{-1}\gamma^{-1} \label{eq:stabhomo} 
\end{align}
is a homomorphism, with kernel $G_{x,\gamma}$. Under the identification $\mathbf{H}\cong T_x$ given by the ordered eigenspaces, the stabilizer of $[(x,\gamma)]$ is $\mathrm{Im}\chi_\gamma$. Indeed, an element identifying $(x,t\gamma)$ with $(x,\gamma)$ must preserve their common pair $(x,y)$. It satisfies $t\gamma = g\gamma g^{-1}$. The action is therefore free exactly when $G_{x,y} = G_{x,\gamma}$.

 It remains to prove the equivalence of (1) and (2). 

Suppose first that \eqref{hyp:spectrum} fails. Choose $I,J,K\subset[n]$, with $|I|=|J|=|K|=r$, $0<r<n$, and $I\ne J$, for which both Horn sets are nonempty. Put $D=\mathrm{diag}\lambda$ and $W=\mathrm{span}\{e_i:i\in I\}$.  In each of the two Horn sets, pick a pair of matrices and conjugate both so that the first becomes $D|_W$ or $D|_{W^\perp}$. The second matrix then defines a Hermitian matrix $B$ on $W\oplus W^\perp$ such that \begin{align*}
\mathrm{eig}(B|_W)&=(-\lambda[J])^\downarrow, \quad\mathrm{eig}(B|_{W^\perp}) = (-\lambda[J^c])^\downarrow,\\
\mathrm{eig}((D+B)|_W)&=\mu[K],\quad \mathrm{eig}((D+B)|_{W^\perp}) =\mu[K^c]. \end{align*}
Since $-B$ and $D$ have the same eigenvalues (counted with multiplicity), there exists a unitary matrix $\gamma$ with  $B =-\gamma^{-1}D\gamma$. After rescaling $\gamma$ by a scalar, we may assume $\det(\gamma)=1$. Define $x:=\ii D$ and $y:=\ii B$. Then \begin{align*}
J(x,\gamma)=\ii(D+B)\in\mathcal{O}_\xi.
\end{align*}

The subspace $W$ reduces both $x$ and $y$. Since $D$ has a simple spectrum, the eigenvalues of $-B|_W$ imply \begin{align*} 
\gamma W=\mathrm{span}\{e_j:j\in J\}\ne W. 
\end{align*}
Let $P_W$ denote the orthogonal projection onto $W$. For $\theta\notin2\pi\mathbb Z$, the matrix \begin{align} g_\theta=\exp\!\left(\ii\theta\left(P_W-\frac rnI_n\right)\right)  \label{eq:block-stabilizer}
\end{align}
belongs to $\mathrm{SU}(n)$ and has the two eigenspaces $W,W^\perp$. It commutes with $x$ and $y$, but not with $\gamma$. Thus, $G_{x,y}\ne G_{x,\gamma}$.

Conversely, suppose that $g\in G_{x,y}\setminus G_{x,\gamma}$. Some eigenspace $W$ of $g$ is not preserved by $\gamma$. Otherwise, $\gamma$ would preserve all eigenspaces and commute with $g$. This subspace is nonzero and proper. Since $g$ commutes with the normal matrices $x,y$, the subspaces $W,W^\perp$ reduce both matrices. Moreover, $\gamma W$ reduces $x$. The simple spectrum of $x$ gives distinct index sets $I,J$, of the same size, such that $W=E_I(x)$ and $\gamma W=E_J(x)$. The restrictions of $-ix,-iy$ to $W$ have spectra $\lambda[I],(-\lambda[J])^\downarrow$. Their sum has spectrum $\mu[K]$ for some set of positions $K$.  The restrictions to $W^\perp$ solve the complementary Horn problem. Repeated entries of $\mu$ are counted with multiplicity. Thus, \eqref{hyp:spectrum} fails. 
\end{proof}

 Although $Z_n\subset\mathbf H$ acts trivially by conjugation, it is not an ineffective subgroup of the spectral action: for $cI_n\in Z_n$, the spectral action sends $[(x,\gamma)]$ to $[(x,c\gamma)]$, and on the free locus, this moves every point unless $c=1$. So the spectral torus is $\mathbf{H}$ itself, with the integral lattice used in Theorem \ref{thm:integral-wall-quotient}.

\begin{remark} \label{rem:optimality}
 By Lemma \ref{lem:hamiltonian-torus} below, the spectral torus acts as a Hamiltonian action on every simple-spectrum orbit-type stratum, even when \eqref{hyp:spectrum} fails. The rank of $h$ at a point $m$ is $n-1-\dim\mathbf{H}_m$, and the stabilizer $\mathbf{H}_m$ is either trivial or of positive dimension. So \eqref{hyp:spectrum} states exactly that $h$ has full rank at every point of every stratum meeting the fiber. Over excluded spectra, chains of lower rank may still exist. The excluded simple spectra lie on the finitely many hyperplanes $\lambda_I-\lambda_J=\mu_K$. However, we exclude only those portions of these hyperplanes where both complementary Horn problems admit solutions. Consequently, the connected components of $\Omega_\xi$ need not coincide with the cells of the resulting hyperplane arrangement. 
\end{remark}

The orbit spectrum $\mu$ is part of the definition of the model, while $\lambda(x)$ is a conserved quantity. Condition \eqref{hyp:spectrum} is therefore a condition on the initial data which is preserved by the flow. It does not change the Hamiltonians or first integrals, and it says nothing about the rational independence of the frequencies of a single trajectory. 

\subsection{Comparison of the two conditions}
\label{subsec:twoconditions}

For two Hermitian $2\times2$ matrices with eigenvalue differences $d_1,d_2$, the possible eigenvalue differences $d_3$ of their sum are exactly
\begin{align}\label{eq:two-dimensional-gap-test} |d_1-d_2|\le d_3\le d_1+d_2. \end{align}
Together with the trace identity, this is the \emph{gap test} used below. Indeed, after subtracting scalar matrices, the squared difference is $d_1^2+d_2^2+2d_1d_2\cos\theta$, with $\cos\theta\in[-1,1]$. 

\begin{proposition}
\label{prop:comparison}
Let $n\ge2$ and fix an orbit spectrum $\mu$. Condition
\eqref{hyp:trace} implies \eqref{hyp:spectrum}. For $n=3$ and
$\mu=(3,-1,-2)$, the feasible spectrum $\lambda=(3,0,-3)$ satisfies
\eqref{hyp:spectrum} and fails \eqref{hyp:trace}.
\end{proposition}

\begin{proof}
 The trace argument in Section \ref{subsec:spectralcond} proves \eqref{hyp:trace}$\Rightarrow$\eqref{hyp:spectrum}. 

 For the converse, we provide an example where the condition (H) holds, but the condition (H') fails.  Let
\begin{align*}
C=\mathrm{diag}(3,-1,-2),\qquad A=\frac14\begin{pmatrix}6&0&9\ii\\0&-2&\ii\sqrt{35}\\-9\ii&-\ii\sqrt{35}&-4\end{pmatrix},\qquad B=A-C=-\overline A.
\end{align*}
 Here $\overline A$ is entrywise conjugation. Direct calculation gives $\mathrm{tr}A=0$, $\mathrm{tr}(A^2)=18$, and $\det A=0$, so the characteristic polynomial is $t^3-9t$. Both $A$ and $B=-\overline A$ therefore have the spectrum $(3,0,-3)$. Choose $\gamma\in G$ with $B=\gamma^*A\gamma$, adjusting a unitary conjurator by a scalar. For $x=\ii A$ and $y=-\ii B$, $J(x,\gamma)=\ii C$. 

In particular, the above construction produces a point of $\mathcal R_\xi$ with variable spectrum $\lambda=(3,0,-3)$ and fixed-orbit spectrum $\mu=(3,-1,-2)$, so $\lambda\in\mathcal B_\xi$ is feasible. However, $\lambda$ does not satisfy \eqref{hyp:trace}: taking $r=1$ and $I=\{1\}$, $J=\{2\}$, $K=\{1\}$ gives \begin{align*}
\lambda_I-\lambda_J=\lambda_1-\lambda_2=3=\mu_1=\mu_K,
\end{align*}
 contradicting the requirement $\lambda_I-\lambda_J\neq\mu_K$ in \eqref{hyp:trace}.

We next check \eqref{hyp:spectrum}. For $n=3$, one only has $r\in\{1,2\}$. The case $r=2$ is equivalent to $r=1$ by passing to complements: the Horn problems for $(I,J,K)$ with $|I|=2$ are exactly the complementary Horn problems for $(I^c,J^c,K^c)$ with $|I^c|=1$ (\eqref{hyp:spectrum}). Thus, it is enough to treat the case $r=1$, i.e., $I=\{i\}$, $J=\{j\}$, $K=\{k\}$. The trace identity is then $\lambda_i-\lambda_j=\mu_k$, with $i\ne j$. Its only solutions are
\begin{align*}
(i,j,k)=(1,2,1)\quad\text{or}\quad(2,3,1).
\end{align*}
In either case, the two Hermitian matrices in the complementary $2\times2$ Horn problem have eigenvalue gaps $3$ and $6$. The required sum has eigenvalues $-1,-2$, so its gap would be $1$.

By the two-dimensional gap test \eqref{eq:two-dimensional-gap-test}, such a sum has a gap of at least $|6-3|=3$, contradicting the required gap $1$. For neither triple are both Horn problems solvable, so \eqref{hyp:spectrum} holds.
\end{proof}

\subsection{Nonemptiness and the Hamiltonian action}
The next two lemmas prove that the admissible domain is dense and the pointwise spectral-rank formula.

\begin{lemma}\label{lem:nonempty-spectrum}
For every $\xi\in\mathfrak{su}(n)$, the set $\Omega_\xi$ is nonempty. The subset $\mathcal{M}_\xi=h^{-1}(\Omega_\xi)$ is open and dense in $\mathcal{R}_\xi$, and $\Omega_\xi$ is dense in $\mathcal{B}_\xi$. In fact, the points satisfying \eqref{hyp:trace} are dense in $\mathcal{R}_\xi$.
\end{lemma}

\begin{proof}

Fix $r\in\{1,\ldots,n-1\}$ and index sets $I,J,K\subset[n]$ with $|I|=|J|=|K|=r$ and $I\neq J$.
There are only finitely many such triples.
By \cite[Theorem 1]{Fulton}, solvability of each (complementary) Horn problem is governed by finitely many linear inequalities (together with the trace constraint) in the spectral variables $\lambda$.
Hence, the set of $\lambda\in\mathcal C^\circ$ for which \emph{both} complementary Horn problems are solvable is closed.
Taking the union over all triples gives a closed subset of $\mathcal C^\circ$, so its complement is open.
Intersecting with the feasible region $\mathcal B_\xi$ (which is itself cut out by finitely many affine inequalities, again by Horn theory) shows that $\Omega_\xi$ is relatively open in $\mathcal B_\xi\cap\mathcal C^\circ$.
Therefore $\mathcal M_\xi=h^{-1}(\Omega_\xi)$ is open in $\mathcal R_\xi$.

First we show that $\mathcal R_\xi$ is nonempty by constructing an element of $J^{-1}(\mathcal O_\xi)$. Let $U\in\mathrm U(n)$ be the standard unitary Fourier matrix with entries
$U_{ab}=n^{-1/2}\exp\bigl(2\pi\ii(a-1)(b-1)/n\bigr)$ \cite[Section 2.2]{HJ}.
For any diagonal $D=\mathrm{diag}(d_1,\ldots,d_n)$, each diagonal entry of $UDU^*$ equals the average $\frac1n\sum_j d_j$.
Since $\sum_j\mu_j=0$, the matrix $Z:=U\mathrm{diag}(\mu)U^*$ has zero diagonal.
Choose distinct unit complex numbers $\varepsilon_j$ with $\prod_j\varepsilon_j=1$, set $\gamma_0=\mathrm{diag}(\varepsilon_1,\ldots,\varepsilon_n)$, and define
\begin{align*}
 C_{ii}=0,\qquad C_{ij}=\frac{Z_{ij}}{1-\overline{\varepsilon_i}\varepsilon_j}\quad(i\ne j).
\end{align*}
Then $C$ is Hermitian, $\mathrm{tr}(C)=0$, and $C-\gamma_0^*C\gamma_0=Z$.
Thus, $(\ii C,\gamma_0)\in J^{-1}(\mathcal O_\xi)$, so $\mathcal R_\xi$ (and hence $\mathcal B_\xi$) is nonempty.

Finally, we show that the points satisfying $(H')$ are dense.
Fix an arbitrary $(x,\gamma)\in J^{-1}(\mathcal O_\xi)$ and set $A:=-\ii x$.
We now vary $A$ within the same moment-map level set: we add a matrix $D_0$ that commutes with $\gamma$, so that the combination $A-\gamma^*A\gamma$ (and hence $J(\ii A,\gamma)$) stays unchanged. Define
\begin{align*}
 d_j:=2^{j-1}-\frac{2^n-1}{n},\qquad j\in[n].
\end{align*}
These numbers sum to zero, and distinct subsets of equal size have distinct sums (by uniqueness of binary expansion).
Choose $u\in\mathrm U(n)$ diagonalizing $\gamma$, and let $D_0$ be the diagonal Hermitian matrix with diagonal entries $d_1,\ldots,d_n$ in that basis.
Conjugating back if needed, we may assume $D_0$ commutes with $\gamma$ and has trace zero.
For $t\in\mathbb R$, set $A_t:=A+tD_0$. As $D_0$ commutes with $\gamma$,
\begin{align*}
A_t-\gamma^*A_t\gamma = A-\gamma^*A\gamma,
\end{align*}
so $J(\ii A_t,\gamma)=J(x,\gamma)$ for all $t$.

For $\eta\in\mathbb R^n$ write $\eta_I=\sum_{i\in I}\eta_i$ and set
\[
Q_\mu(\eta)=
\prod_{p<q}(\eta_p-\eta_q)^2
\prod_{r=1}^{n-1}
\prod_{\substack{I,J,K\subset[n]\\|I|=|J|=|K|=r,\ I\ne J}}
(\eta_I-\eta_J-\mu_K).
\]
Permuting the variables permutes the factors, so $Q_\mu$ is symmetric.
Every symmetric polynomial in the eigenvalues is a polynomial in the
coefficients of the characteristic polynomial. Consequently
\[
P(t):=Q_\mu(\mathrm{eig}(A+tD_0))
\]
is a polynomial in $t$. Put $N=\deg Q_\mu$. Its highest-degree
homogeneous part is $Q_0$, and continuity of the ordered eigenvalues,
together with symmetry, gives
\[
\lim_{t\to+\infty}t^{-N}P(t)=Q_0(d_1,\ldots,d_n)\ne0.
\]
The last inequality follows because the $d_j$ are distinct and
distinct subsets of equal cardinality have different sums. Thus
$P$ is not identically zero and has only finitely many zeros.
Choose $t_\nu\to0$ with $P(t_\nu)\ne0$. The spectrum of $A+t_\nu D_0$
is then simple and satisfies \eqref{hyp:trace}, and
\[
J(\ii(A+t_\nu D_0),\gamma)=J(x,\gamma),\qquad
(\ii(A+t_\nu D_0),\gamma)\longrightarrow(x,\gamma).
\]
Passing to the quotient proves density in $\mathcal R_\xi$. Since
\eqref{hyp:trace} implies \eqref{hyp:spectrum} and the spectral map
is continuous and surjective onto $\mathcal B_\xi$, this proves the
remaining nonemptiness and density assertions.

\end{proof}

\begin{lemma}\label{lem:hamiltonian-torus}
Let $\mathbf{H}$ act on the simple-spectrum part of $\mathcal R_\xi$ by \eqref{eq:eigenvalue-action}, and let $S$ be a connected component of the intersection of an original symplectic orbit-type stratum with $h^{-1}(\mathcal{C}^\circ)$. Then the action preserves $S$ and restricts to a Hamiltonian action on $S$ with moment map $h|_S$. For $m=[(x,\gamma)]\in S$ with stabilizer $\mathbf{H}_m$,
\begin{align}\label{eq:pointwisespec}
\dim\mathbf{H}_m=\dim G_{x,y}-\dim G_{x,\gamma},\qquad
\mathrm{rank}d_m(h|_S)=n-1-\dim\mathbf{H}_m,
\end{align}
and $\mathbf{H}_m$ is either trivial or positive dimensional. On $\mathcal{M}_\xi$ the action is free, $\mathrm{rank}d(h|_S)=k=n-1$ on every stratum $S\subset\mathcal{M}_\xi$, and $\Omega_\xi$ is open in $\mathcal{C}^\circ$. 
\end{lemma}
\begin{proof}
Fix $a\in\R^n$ with $\sum_i a_i=0$ and consider the spectral Hamiltonian $h_a=\sum_i a_i\lambda_i$. Its Hamiltonian flow on $T^*G\cong\mathfrak{g}\times G$ is given by $\dot x=0$ and $\dot\gamma=\ii\sum_i a_iE_i(x)\gamma$ \cite[Lemma 3.6]{FF2026}, i.e. the action \eqref{eq:eigenvalue-action}. In particular it preserves $x$ (hence the orbit type of $x$), preserves $J$ (so it descends to the shifted reduction), and therefore preserves the intersection of an orbit-type stratum with $h^{-1}(\mathcal{C}^\circ)$. Since $G_x$ is abelian, the stabilizer $G_{x,\gamma}$ is unchanged under multiplication by $t(x)\in\mathbf H$, so the induced $\mathbf H$-action preserves each such connected component $S$. By construction, the corresponding moment map is $h|_S$.

For $m=[(x,\gamma)]\in S$, the stabilizer dimension formula in \eqref{eq:pointwisespec} follows from the homomorphism \eqref{eq:stabhomo} by comparing its kernel and image. Moreover, for a Hamiltonian torus action one has $\mathrm{rank}\,d_m(h|_S)=\dim(\mathbf H\cdot m)= (n-1)-\dim\mathbf H_m$, which is the second identity.

Finally, if $\mathbf H_m$ is nontrivial, pick $g\in G_{x,y}\setminus G_{x,\gamma}$. As in the proof of Proposition \ref{prop:stabilizer}, choose an eigenspace $W$ of $g$ not preserved by $\gamma$ and the associated circle $g_\theta$ from \eqref{eq:block-stabilizer}. Then
\begin{align*}
\chi_\gamma(g_\theta)=\exp\!\left(\ii\theta(P_W-P_{\gamma W})\right)
\end{align*}
is nonconstant, so $\mathbf H_m$ contains a circle and is therefore positive dimensional.

On $\mathcal M_\xi$, Proposition \ref{prop:stabilizer} shows that the $\mathbf H$-action is free, hence $\mathrm{rank}\,d(h|_S)=n-1$ on every stratum $S\subset\mathcal M_\xi$. The openness of $\Omega_\xi\subset\mathcal{C}^\circ$ then follows by applying the submersion theorem on a stratum through a preimage of each point of $\Omega_\xi$.
\end{proof}

\subsection{The quotient topology, smooth functions, and strata}
The stabilizer calculation identifies the fibers of $\psi$ as
spectral-torus orbits. To use this identification in the stratified
Poisson chain, we must also identify the quotient topology, smooth
functions, and orbit-type strata.

\begin{lemma}\label{lem:global-torus-quotient}
Let $\mathbf{H}$ act on $\mathcal{M}_\xi$ by \eqref{eq:eigenvalue-action}, and let $\psi:\mathcal{M}_\xi\to\mathcal{P}_\xi^\Omega$ be the map defined above. Suppose $\mathcal{M}_\xi$ and $\mathcal{P}_\xi^\Omega$ are equipped with their quotient topologies and their quotient smooth structures (Section \ref{subsec:orbit-types}).
Then:
\begin{enumerate}
\item The map $\psi$ is $\mathbf{H}$-invariant and induces a map $\bar\psi:\mathcal{M}_\xi/\mathbf{H}\to\mathcal{P}_\xi^\Omega$. Moreover, $\bar\psi$ is a homeomorphism, and it identifies the smooth functions on the two spaces.

\item Let $p\in\mathcal{P}_\xi^\Omega$.
Then there exist an open neighborhood $U\subset\mathcal{P}_\xi^\Omega$ of $p$ and an $\mathbf{H}$-equivariant isomorphism of smooth stratified spaces \begin{align}
 \Phi:\psi^{-1}(U)\xrightarrow{\ \cong\ }\mathbf{H}\times U, \label{eq:torlocpro}
\end{align}
where $\mathbf{H}$ acts on $\mathbf{H}\times U$ by left multiplication on the first factor, and such that $\psi=\mathrm{pr}_2\circ\Phi$. This holds also when $p$ lies in a singular stratum.

\item Let $\{S_\alpha\}$ be the connected symplectic orbit-type strata of $\mathcal{M}_\xi$.
Then the collection $\{\psi(S_\alpha)\}$ is a finite Whitney stratification of $\mathcal{P}_\xi^\Omega$.
\end{enumerate}
\end{lemma}

\begin{proof}
Put $N_0=\mathfrak{g}_{\rm reg}\times G$ and $E=\{(x,y):x\in\mathfrak{g}_{\rm reg},\ y\in\mathcal{O}_{-x}\}$. Ordered local eigenframes embed $E$ smoothly in $\mathfrak{g} \times \mathfrak{g}$ and give local sections of the principal $\mathbf{H}$-bundle $F:N_0\to E$, $F(x,\gamma)=(x,-\gamma^{-1}x\gamma)$. Let $D\subset E$ impose $x+y\in\mathcal{O}_\xi$ and $\lambda(x)\in\Omega_\xi$, and set $S=F^{-1}(D)$. The commuting compact actions give
\begin{equation}\label{eq:quotient-identification}
\mathcal{M}_\xi/\mathbf{H}=(S/G)/\mathbf{H}\cong S/(G\times\mathbf{H})
\cong (S/\mathbf{H})/G=D/G=\mathcal{P}_\xi^\Omega
\end{equation}
with the quotient topologies. The principal local sections identify $S/\mathbf{H}$ with $D$ in its subspace topology.

Pullback by $F$ takes smooth functions on $D/G$ to $\mathbf{H}$-invariant smooth functions on $S/G$. Conversely, patch local extensions of such a function near a compact $(G\times\mathbf{H})$-orbit and average over this group. The resulting invariant smooth function on $N_0$ descends through $F$ to a smooth function on $E$, and then extends locally to $\mathfrak{g} \times \mathfrak{g}$. Averaging over $G$ gives the required invariant extension. Thus \eqref{eq:quotient-identification} also identifies smooth structures.

Fix a point $m\in S\subset N_0$ and write $L:=G_m$ for its stabilizer under the $G$-action. Since the $\mathbf H$-action on $S/G$ is free, the stabilizer of $m$ for the product action of $G\times\mathbf H$ is $(G\times\mathbf H)_m=L\times\{1\}$. By the slice theorem for compact group actions \cite[Theorem 1.25]{Meinrenken}, the $G\times\mathbf H$-orbit through $m$ has a standard local model. More precisely, there is a $G\times\mathbf H$-invariant open neighborhood $\mathcal U\subset N_0$ of $(G\times\mathbf H)\cdot m$ and a $G\times\mathbf H$-equivariant diffeomorphism
\begin{align*}
\mathcal U \; \longrightarrow\;  (G\times\mathbf H)\times_L B,
\end{align*}
where $B$ is a small $L$-invariant ball in the normal space $\nu_m:=T_mN_0/T_m((G\times\mathbf H)\cdot m)$.
Here $L$ acts on $G\times\mathbf H$ by right multiplication and on $B$ by the induced linear (slice) representation.

Quotienting this local model by $G$ simply forgets the $G$-coordinate.
This gives a local description of $N_0/G$ near the class $[m]$ as a product of the torus factor and the slice quotient:
\begin{equation}\label{eq:slice-product}
((G\times\mathbf{H})\times_L B)/G\cong\mathbf{H}\times(B/L),
\qquad [g,t,v]\longmapsto\bigl(t,[v]\bigr).
\end{equation}

The subset $S\subset N_0$ is cut out by the constraints defining $D\subset E$.
In the slice model, this means that, after shrinking $B$ if needed, $S$ corresponds to $(G\times\mathbf H)\times_L C_0$ for some $L$-invariant subset $C_0\subset B$. Therefore, in the same local coordinates,
\begin{align*}
S/G &\cong \mathbf H\times (C_0/L).
\end{align*}
With the quotient smooth structures (defined via invariant functions), this identification is smooth. 
Finally, dividing by $\mathbf H$ identifies a neighborhood in $S/(G\times\mathbf H)$ with an open subset $U\subset D/G$, which gives \eqref{eq:torlocpro}.

In these charts orbit-type strata are products $\mathbf{H}\times A_\tau$. The quotient smooth structure agrees with the reduction charts of \cite[Theorems 7 and 15]{BL}. Their transitions have local smooth ambient extensions in both directions. Whitney conditions are unchanged by these coordinate changes or by taking a transverse slice to the smooth torus factor. Thus the $A_\tau$ are Whitney strata. Connectedness of $\mathbf{H}$ makes each orbit-type stratum of $\mathcal{M}_\xi$ a union of $\psi$-fibers, so distinct strata have disjoint images. A set $T$ is $\psi$-saturated if $\psi^{-1}(\psi(T))=T$. Since $\psi$ is open, for such a set $T$,
$\psi^{-1}(\overline{\psi(T)})=\overline T$.
This transfers the frontier condition to $\mathcal{P}_\xi^\Omega$ and proves assertion (3), apart from finiteness.

For finiteness, the unital $*$-algebra generated by $x,\gamma$ decomposes its representation as
$\mathbb C^n=\bigoplus_j(\mathbb C^{d_j}\otimes\mathbb C^{m_j})$, with $\sum_jd_jm_j=n$.
Its unitary commutant is $\prod_j(I_{d_j}\otimes\mathrm U(m_j))$, and intersecting with $G$ imposes $\prod_j\det(u_j)^{d_j}=1$. There are finitely many possible integer lists, hence stabilizer types. For each type, the condition that the stabilizer be conjugate to it is a first-order condition on real matrix entries, so is semialgebraic by quantifier elimination. The constraints and $\Omega_\xi$ are semialgebraic too. Their images in invariant coordinates consequently have finitely many connected components \cite[Sections 2.2 and 2.4]{BCR}. The same argument applies to the pairs $(x,y)$.
\end{proof}

Write $\mathcal{M}_\alpha$ for these connected orbit-type strata of $\mathcal{M}_\xi$ and $\mathcal{P}_\delta$ for their images. In particular,
\begin{align}\label{eq:matched-strata} \mathcal{P}_\delta=\mathcal{M}_\alpha/\mathbf{H},\qquad \psi^{-1}(\mathcal{P}_\delta)=\mathcal{M}_\alpha. \end{align}
Let $\mathcal{B}_\beta$ be the connected components of $\Omega_\xi$. They are finitely many open subsets of the trace-zero hyperplane. Every $\mathcal{P}_\delta$ maps into one $\mathcal{B}_\beta$.

\subsection{Poisson structures and the subchains}
We now verify the Poisson and rank conditions.

\begin{lemma}\label{lem:poisson-ranks}
The natural brackets on $\mathcal{M}_\xi$ and $\mathcal{P}_\xi^\Omega$ are represented by smooth stratified Poisson bivectors tangent to the strata above. The maps $\psi,\pi$ are Poisson, and the base bracket is zero. If $d_\alpha=\dim\mathcal{M}_\alpha$, then
\begin{align}\label{eq:subchain-ranks} \dim\mathcal{P}_\delta=d_\alpha-k,\qquad \mathrm{rank}\Pi_\delta=d_\alpha-2k,\qquad \mathrm{rank}(d\pi|_{\mathcal{P}_\delta})=k. \end{align}
The chain $\mathcal{M}_\alpha\to\mathcal{P}_\delta\to\mathcal{B}_\beta$ is superintegrable in the sense of Definition \ref{def:strong}.
\end{lemma}
\begin{proof}
First, the bracket on $\mathcal{M}_\xi$ is the reduced bracket \eqref{eq:reduced-bracket} coming from singular reduction, and the bracket on $\mathcal{P}_\xi^\Omega$ is the quotient bracket for the Hamiltonian $\mathbf{H}$-action on $\mathcal{M}_\xi$. This quotient bracket agrees with the pair bracket \eqref{eq:product-poisson}, since the map $(x,\gamma)\mapsto (x,-\gamma^{-1}x\gamma)$ is Poisson from $T^*G$ to its two momentum factors; in particular $\psi$ is Poisson.

Next, choose real polynomial generators $p_1,\ldots,p_N$ for the invariant algebra on $\mathfrak{g}\times\mathfrak{g}$. Their Poisson brackets have the form $\{p_i,p_j\}=b_{ij}(p_1,\ldots,p_N)$ for polynomials $b_{ij}$, so the bivector $\sum_{i<j} b_{ij}(u)\,\partial_{u_i}\wedge\partial_{u_j}$ represents the quotient bracket in invariant coordinates. Invariant Hamiltonian flows preserve both spectra, $x+y$, and simultaneous stabilizers; hence each column of $(b_{ij})$ is tangent to every quotient stratum, and by skew symmetry the bivector is tangent as well. Schwarz's theorem \cite{Schwarz} extends this representation from polynomial invariants to all smooth invariant functions, and Jacobi follows from the Jacobi identity on $\mathfrak g\times\mathfrak g$. The same argument applies to $T^*G/G$: the cotangent bivector has polynomial coefficients in $(x,\gamma)$, averaging produces invariant polynomial extensions, and invariant Hamiltonian flows fix $J$ and preserve stabilizers. Therefore the resulting bivector restricts to $\mathcal{M}_\xi$ and is tangent to its orbit-type strata.

We then show that $\pi$ is Poisson with zero bracket on the base. Indeed, the Hamiltonian flows of the eigenvalue functions $\lambda_i$ generate the spectral-torus action, so each $\lambda_i$ Poisson-commutes with every $\mathbf{H}$-invariant function.

Finally, we prove the rank identities \eqref{eq:subchain-ranks} and identify the symplectic leaves. Fix $m\in\mathcal{M}_\alpha$ and set $V_m=T_m(\mathbf{H}m)$. For a subspace $V$ of a symplectic vector space, write $V^\omega=\{w:\omega(w,v)=0\ \forall v\in V\}$, and for a bivector $\Pi$, define $\Pi^\sharp$ by $\beta(\Pi^\sharp\alpha)=\Pi(\alpha,\beta)$. Since the $\mathbf{H}$-action on $\mathcal{M}_\alpha$ is free and $\mathbf{H}$ is abelian, the moment-map identity gives
\begin{align}\label{eq:vertical-orthogonal}
\dim V_m=k,\qquad
V_m\subset V_m^{\omega_\alpha},\qquad
\ker d_m(h|_{\mathcal{M}_\alpha})=V_m^{\omega_\alpha}.
\end{align}
The quotient submersion $\psi_\alpha:\mathcal{M}_\alpha\to\mathcal{P}_\delta$ has kernel $V_m$, and since $h=\pi\circ\psi$, it follows that $\mathrm{rank}(d\pi|_{\mathcal{P}_\delta})=k$. Pullback by $\psi_\alpha$ identifies covectors on $\mathcal{P}_\delta$ with the annihilator of $V_m$, and the Poisson identity yields
\begin{align}\label{eq:poisson-image}
\mathrm{Im}\Pi_\delta^\sharp
=d\psi_\alpha(V_m^{\omega_\alpha})
=\ker d_{\psi(m)}\pi.
\end{align}
Since $\dim V_m^{\omega_\alpha}=d_\alpha-k$, quotienting by $V_m$ gives the remaining identities in \eqref{eq:subchain-ranks}. Equation \eqref{eq:poisson-image} also identifies the symplectic foliation of $\mathcal{P}_\delta$ with the connected components of the fibers of the constant-rank map $\pi|_{\mathcal{P}_\delta}$. Therefore, the chain $\mathcal{M}_\alpha\to\mathcal{P}_\delta\to\mathcal{B}_\beta$ is superintegrable in the sense of Definition~\ref{def:strong}.
\end{proof}

\begin{proof}[Proof of Theorem \ref{thm:main}]

We proceed in four steps corresponding to parts \textnormal{(i)}-\textnormal{(iv)}. First we establish the description of the free locus and the required Whitney stratifications. Next we verify that the resulting maps form a stratified Poisson chain and identify the (trivial) characteristic quotient on this locus. We then prove the principal-bundle and properness statements. Finally, we check completeness of the commuting Hamiltonian flows.

\smallskip
\noindent\emph{Part \textnormal{(i)} (free locus and stratifications).}
First, the two equivalences are precisely Proposition \ref{prop:stabilizer}.

Next, Lemma \ref{lem:nonempty-spectrum} shows that $\Omega_\xi\subset\mathcal{B}_\xi$ is nonempty and dense, and that $\mathcal{M}_\xi=h^{-1}(\Omega_\xi)$ is open and dense in $\mathcal{R}_\xi$. Moreover, Lemma \ref{lem:hamiltonian-torus} implies that $\Omega_\xi$ is open in $\mathcal{C}^\circ$.

Finally, we show that $\Omega_\xi$ is semialgebraic with a finite Whitney stratification. By the Horn theorem, each Horn set in \eqref{hyp:spectrum} is nonempty exactly when finitely many affine inequalities in $\lambda$ hold, and $\mathcal{B}_\xi$ is a convex polyhedron (Section \ref{subsec:observable}). Hence, $\Omega_\xi$ is semialgebraic and has finitely many connected components \cite[Section 2.4]{BCR}; each component is open in the trace-zero hyperplane and therefore yields a finite Whitney stratification of $\Omega_\xi$. Lemma \ref{lem:global-torus-quotient} (and its proof) then supplies the finite Whitney stratifications of $\mathcal{M}_\xi$ by connected orbit-type strata and of $\mathcal{P}_\xi^\Omega$ by their images.

\smallskip
\noindent\emph{Part \textnormal{(ii)} (Poisson chain and ``no refinement'' case).}
Lemma \ref{lem:poisson-ranks} shows that $\psi$ and $\pi$ are Poisson, that the bracket on the base is zero, and that the stratal symplectic forms are Poisson compatible (because the stratified bivectors are tangent to the strata and restrict to the symplectic brackets).
Surjectivity of $\psi$ is immediate from the definition of $\mathcal{P}_\xi^\Omega$.
This gives the original chain of Definition \ref{def:stratified-system}.

For the compatible refinements and global characteristic quotient in Definition \ref{def:stratified-system}, nothing changes here: each middle stratum $\mathcal{P}_\delta$ is connected and maps into a single component $\mathcal{B}_\beta$ of $\Omega_\xi$, so $T=\mathcal{P}_\delta$.
Then \eqref{eq:matched-strata} gives $\mathcal{M}_\alpha\cap\psi^{-1}(T)=\mathcal{M}_\alpha$, which is connected, hence $S=\mathcal{M}_\alpha$.
Therefore $\tau_S=\omega_\alpha$ is symplectic, so its characteristic leaves are points and the relation \eqref{eq:separate-relation} is just equality.
Consequently $p=\mathrm{id}$, $\underline{\mathcal{M}}=\mathcal{M}_\xi$ (with its own topology and smooth structure), and $\underline\psi=\psi$.
The diagrams \eqref{eq:global-reduction-diagram} and \eqref{eq:stratum-reduction-diagram} commute, and their bottom rows are strong by Lemma \ref{lem:poisson-ranks}.
This proves that \eqref{eq:main-diagram} is stratified superintegrable.
The final closure statement follows because, for $S\subset\mathcal{M}_\xi$, the closure in the subspace topology is $\mathcal{M}_\xi\cap\overline{S}^{\,\mathcal{R}_\xi}$.

\smallskip
\noindent\emph{Part \textnormal{(iii)} (principal bundle and properness).}
Lemma \ref{lem:hamiltonian-torus} gives that the $\mathbf{H}$-action on $\mathcal{M}_\xi$ is free and Hamiltonian with moment map $h$.
By Lemma \ref{lem:global-torus-quotient}(1)-(2), $\psi$ is the quotient map for this action and is locally a product with $\mathbf{H}$, including at singular strata.
The map $\pi$ is surjective because every feasible spectrum occurs on $\mathcal{R}_\xi$, and points with spectrum in $\Omega_\xi$ lie in $\mathcal{M}_\xi=h^{-1}(\Omega_\xi)$.

For properness, let $K\subset\Omega_\xi$ be compact.
Before quotienting, $h^{-1}(K)$ is closed and compact since $G$ and $\mathcal{O}_\xi$ are compact and $\|x\|^2=\sum_i\lambda_i(x)^2$ is bounded on $K$.
Hence $h^{-1}(K)$ and $\pi^{-1}(K)=\psi(h^{-1}(K))$ are compact.
If $C\subset\mathcal{P}_\xi^\Omega$ is compact, then $\psi^{-1}(C)$ is closed in $h^{-1}(\pi(C))$, which is compact, so $\psi^{-1}(C)$ is compact.
Thus both maps are proper (without asserting properness on individual strata).

Finally, \eqref{eq:main-dimensions} is \eqref{eq:subchain-ranks}, and Lemma \ref{lem:poisson-ranks} implies that every induced chain is strong and that the connected components of the fibers of $\pi|_{\mathcal{P}_\delta}$ are its symplectic leaves.

\smallskip
\noindent\emph{Part \textnormal{(iv)} (completeness of the commuting flows).}
Equation \eqref{eq:original-flow} proves completeness and shows that the flow preserves $x,y,J$ and the orbit type (hence the spectrum) on every stratum, and preserves $\mathcal{M}_\xi$.
For a smooth function $F(H_2,\ldots,H_n)$, replace $D_j(x)$ by $\sum_{j=2}^n(\partial_jF)(H_2(x),\ldots,H_n(x))D_j(x)$; this is constant along the flow and gives the same conclusion.
\end{proof}

\section{The global characteristic quotient}\label{sec:intewalquo}

We now retain the simple-spectrum points at which the spectral
torus acts freely, without requiring their entire spectral fibers
to be free. Each connected orbit-type stratum of
$\mathcal M_\xi^{\mathrm{free}}$ is symplectic. The selected walls
subdivide it into open sign regions and submanifolds lying on one
or more walls. The form remains symplectic on an open sign region.
On a piece lying on $c$ independent walls, its kernel has dimension
$c$ and is generated by the corresponding characteristic subtorus.
The global characteristic quotient identifies the points of each
such subtorus orbit, while keeping the adjacent regions and all
wall intersections in one space.

\subsection{The free simple-spectrum domain}\label{subsec:pointwise-free}

By Proposition \ref{prop:stabilizer}, the free sets are
\begin{align}\label{eq:pointwise-free}
 \mathcal{M}_\xi^{\mathrm{free}}&=\{[(x,\gamma)]\in\mathcal{R}_\xi:\lambda(x)\in\mathcal{C}^\circ,\ G_{x,y}=G_{x,\gamma}\},\\
 \mathcal{P}_\xi^{\mathrm{free}}&=\psi(\mathcal{M}_\xi^{\mathrm{free}}),\qquad
 \mathcal{B}_\xi^{\mathrm{free}}=h(\mathcal{M}_\xi^{\mathrm{free}}).\nonumber
\end{align}
The same proposition gives
\begin{align*}
 \Omega_\xi=\{\lambda\in\mathcal{B}_\xi\cap\mathcal{C}^\circ:h^{-1}(\lambda)\subset\mathcal{M}_\xi^{\mathrm{free}}\}.
\end{align*}

\begin{proposition}[The free simple-spectrum restriction]\label{prop:pointwise-free}
The sets in \eqref{eq:pointwise-free} are open semialgebraic subsets of the original simple-spectrum spaces, and $\mathcal{M}_\xi^{\mathrm{free}}$ is dense in $\mathcal{R}_\xi$. With the original smooth structures, the chain
\begin{align*}
 \mathcal{M}_\xi^{\mathrm{free}}\xrightarrow{\psi}
 \mathcal{P}_\xi^{\mathrm{free}}\xrightarrow{\pi}
 \mathcal{B}_\xi^{\mathrm{free}}
\end{align*}
is stratified superintegrable in the sense of Definition \ref{def:stratified-system}, with $p=\mathrm{id}$
and rank $n-1$ on every subchain. Use the connected orbit-type strata of $\mathcal{M}_\xi^{\mathrm{free}}$, their images under $\psi$ in $\mathcal{P}_\xi^{\mathrm{free}}$, and the connected components of $\mathcal{B}_\xi^{\mathrm{free}}$ as the strata. The map $\psi$ is a proper principal $\mathbf{H}$-bundle in the local-product sense of Lemma \ref{lem:global-torus-quotient}. The Hamiltonians of Theorem \ref{thm:main}(iv) have complete stratum-preserving flows. The spectral projection $\pi$ need not be proper.
\end{proposition}
\begin{proof}

Proposition \ref{prop:stabilizer} identifies the fibers of $\psi$ on the simple-spectrum domain with $\mathbf{H}$-orbits, so the free locus is saturated. At a free reduced point represented by $m$, the $G\times\mathbf{H}$-stabilizer is $G_m\times\{1\}$. In its compact-group slice, nearby stabilizers have trivial projection to $\mathbf{H}$, proving openness and the local torus products. The quotient, averaging, and semialgebraic arguments of Lemma \ref{lem:global-torus-quotient} then give the restriction smooth structure and finite Whitney stratifications on this saturated open set. Density follows from $\mathcal{M}_\xi\subseteq\mathcal{M}_\xi^{\mathrm{free}}$.

Lemma \ref{lem:hamiltonian-torus} makes $h$ a submersion on each stratum, so $\mathcal{B}_\xi^{\mathrm{free}}$ is open. Lemma \ref{lem:poisson-ranks} gives the strong chains, with zero characteristic distributions and $p=\mathrm{id}$. The quotient map $\psi$ by the compact torus is proper, and the flows \eqref{eq:original-flow} preserve this domain. Remark \ref{rem:free-nonproper} shows that $\pi$ need not be proper.

\end{proof}

\subsection{The local symplectic model}

The following application of the Hamiltonian normal form supplies the local products needed for the simultaneous wall quotient. The global argument will retain both sides of every wall.

\begin{lemma}[A symplectic product near a free spectral orbit]
\label{lem:free-spectral-normal-form}
Every $m\in\mathcal{M}_\xi^{\mathrm{free}}$ has an
$\mathbf{H}$-invariant neighborhood which is isomorphic, as a
stratified space with its quotient smooth structure and its symplectic
orbit-type strata, to a product
\begin{align*}
\mathbf{H}\times U\times Y,\qquad h(\theta,\lambda,y)&=\lambda,\qquad \omega =\omega_{T^*\mathbf{H}}+\omega_Y\ \text{on every stratum}.
\end{align*}
Here $U\subset\mathfrak h^*$ is open. $V$ is the symplectic slice representation of the stabilizer $L$ of a lift of $m$, with quadratic
moment map $\langle\Phi_V(v),X\rangle=\tfrac12\omega_V(v,Xv)$ for $X\in\mathfrak l=\mathrm{Lie}(L)$. $V/\!/L$ denotes $\Phi_V^{-1}(0)/L$. The space $Y$ is a
neighborhood of the image of $0$ in $V/\!/L$, with reduced stratal forms $\omega_Y$. In these coordinates
$\mathcal{P}_\xi^{\mathrm{free}}$ is locally $U\times Y$, and $\psi$
forgets $\theta$.
\end{lemma}

\begin{proof}

On the simple-spectrum part of the shifted manifold $\mathscr M$ of \eqref{eq:shift}, the group $\Gamma=G\times\mathbf{H}$ has a moment map $(\Phi,h-h(m))$. At a lift $a$ of $m$, its value is zero and $\Gamma_a=L\times\{1\}$. Thus $\Gamma a$ is isotropic, with symplectic slice $V=(T_a(\Gamma a))^{\omega_{\mathscr M}}/T_a(\Gamma a)$.

The Hamiltonian normal form \cite[Theorem 9 and Corollary 10]{BL} gives a neighborhood of the zero section in
\begin{align*}
 (G\times\mathbf{H})\times_L
 (\mathrm{ann}(\mathfrak l)\oplus\mathfrak h^*\oplus V).
\end{align*}
where $\mathrm{ann}(\mathfrak l)=\{\eta\in\mathfrak{g}^*: \eta|_{\mathfrak l}=0\}$. An $L$-invariant splitting $\mathfrak{g}^*=\mathrm{ann}(\mathfrak l)\oplus\mathfrak l^*$ identifies the moment map with $(\mathrm{Ad}_g^*(\nu+\Phi_V(v)),s)$, where $s=h-h(m)$. The $G$-moment equation is $\nu=0$, $\Phi_V(v)=0$. Reduction by $G$ therefore gives $\mathbf{H}\times\mathfrak h^*\times(V/\!/L)$. The normal form arises by $L$-reduction of $T^*(G\times\mathbf{H})\times V$, leaving the $T^*\mathbf{H}$ factor unchanged, so the form splits as stated. Reduction also identifies the quotient smooth structures and orbit-type strata. Shrinking gives the required product neighborhood.

\end{proof}

\subsection{The global quotient}
 We now prove Theorem \ref{thm:integral-wall-quotient}, in its
notation. Hausdorffness is established by the closed-saturation argument in Step~2 below.
Proposition \ref{prop:wall-local-model} describes the refined strata
and their characteristic leaves, gives the local model of the
quotient, and proves the dimension and rank formulas. The proof of the
theorem, at the end of this section, combines this argument with the
local models in five steps.

In the local product of Lemma \ref{lem:free-spectral-normal-form},
use moment coordinates $\rho_j=\langle\lambda,b_j\rangle$ and their
dual circle coordinates $\theta_j$. On an original stratum of
dimension $d$,
\[
\omega=\sum_{j=1}^k d\rho_j\wedge d\theta_j+\omega_Y.
\]
On a piece with $c$ active walls, the corresponding $d\rho_j$ vanish
on tangent vectors. The restricted kernel is spanned by those
$\partial_{\theta_j}$, so
\[
\dim S=d-c,\qquad \dim p(S)=d-2c,\qquad
\text{rank of the induced strong chain}=k-c.
\]
Integral normals give periodic Hamiltonian flows and closed
characteristic subtori. A common lattice basis supplies independent
primitive circle factors and the product-of-cones local models
below. This is a sufficient hypothesis for our construction; no
necessity assertion is made for general Hausdorff wall quotients.
Closed individual leaves alone do not establish Hausdorffness;
the closed-saturation argument below supplies the global verification.

The local quotient uses the following double-cone calculation. Give $N=\mathbb R_s\times S^1_\theta$, with $\omega=ds\wedge d\theta$, the strata $N_\pm=\{\pm s>0\}$ and $F=\{s=0\}$. Their characteristic leaves are points on $N_\pm$ and the circle $F$. The map
\begin{equation}
p_C(s,\theta)=(s\cos\theta,s\sin\theta,s),\qquad C=\{(x,y,t)\in\mathbb R^3:x^2+y^2=t^2\}, \label{eq:cone-quotient}
\end{equation}
is proper, since $\|p_C(s,\theta)\|^2=2s^2$, and its fibers are exactly these leaves. Thus $C$ realizes the characteristic quotient. Equip it with the ambient smooth functions and strata $C_\pm=\{\pm t>0\}$ and $\{0\}$. This stratification is Whitney: the tangent planes on $C_\pm$ contain the radial secants to the vertex, and condition (A) at the point stratum is automatic. The reduced forms and coordinate brackets are
\begin{equation}\label{eq:cone-symplectic-form}
\underline\omega_C|_{C_\pm}=\frac{dx\wedge dy}{t},\qquad
\{x,y\}=t,\quad \{t,x\}=-y,\quad \{t,y\}=x.
\end{equation} Indeed $p_C^*(dx\wedge dy)=s\,ds\wedge d\theta$.
The displayed linear brackets satisfy Jacobi, have $x^2+y^2-t^2$ as a Casimir, and vanish at the vertex. They therefore define a
stratified Poisson structure whose restrictions invert the reduced forms. This is the Poisson double cone of \cite[Section 4]{HCES2024}.

\begin{proposition}[Local models for the wall quotient]
\label{prop:wall-local-model}
In the setting of Theorem \ref{thm:integral-wall-quotient}, at a point
on $c=|I|$ walls a local model of $X$ is an open subset of
\begin{equation}\label{eq:wall-cone-local-model}
C^c\times T^*(\mathbf{H}/K_I)\times Y,
\end{equation}
where $C$ is the double cone \eqref{eq:cone-quotient} with
bracket \eqref{eq:cone-symplectic-form}, and $Y$ is the symplectic slice quotient of
Lemma \ref{lem:free-spectral-normal-form}.
On a refined stratum $S$ contained in an original symplectic stratum
of dimension $d$, the characteristic fibers have dimension $c$ and
\begin{align}\label{eq:wall-rank-formulas}
\dim S&=d-c,& \dim p(S)&=d-2c,&
\dim\psi(S)&=d-k-c,\\
\mathrm{rank}\Pi_{\psi(S)}&=d-2k,&
\mathrm{rank}d(\pi|_{\psi(S)})&=k-c.\nonumber
\end{align}
Thus $p$ is not injective on any nonempty refined stratum lying on
at least one wall, and the chain induced on $p(S)$ in
Theorem \ref{thm:integral-wall-quotient} has rank $k-c$.
\end{proposition}

\begin{proof}
\emph{Refined strata.}
Lemma \ref{lem:free-spectral-normal-form} makes $h$ a projection
near every point, including points where several original strata
meet. Each sign refinement is locally
a product of a coordinate hyperplane arrangement with the slice
stratification, both on $M$ and on $P$. These product models preserve
the Whitney conditions.

For the frontier condition, let $S,S'$ be connected refined strata and $A=S'\cap\overline S$. This is closed in $S'$. At $m\in A$, decompose the slice as $V=V^L\oplus(V^L)^\omega$, so $V/\!/L=V^L\times((V^L)^\omega/\!/L)$. Together with the coordinate walls this gives a local product $U_0\times Z$ in which $S'$ is $U_0\times\{0\}$ and $U_0$ is connected. The remaining local branches are $U_0\times Z_j$ with finitely many connected $Z_j$, by semialgebraicity. A branch meeting $S$ lies in $S$. Since $m\in\overline S$, one has $0\in\overline{Z_j}$ for such a branch. Hence $A$ contains $U_0\times\{0\}$ and is open in $S'$. Connectedness gives $A=S'$. Omitting the torus factor proves the middle frontier condition. The same argument applies to connected sign pieces in the base.

The refinement has finitely many strata, since the original orbit-type strata, the domain and the wall equations are all
semialgebraic. 

\emph{Characteristic leaves.}
On a refined stratum $S$ with active set $I$, the moment map identity gives  \begin{align*}
 \ker(\omega|_S)  &= \mathrm{span}\{(b_j)_M:j\in I\},  \quad  (b_j)_M(m) =\left.\frac{d}{ds}\right|_{s=0}\exp(\ii s\mathrm{diag}b_j)\cdot m.
\end{align*}
Since the action is free, this distribution has dimension $c$, so $\omega|_S$ is presymplectic. By
the assumption on the basis, $K_I$ is a closed connected subtorus,
and its orbits are exactly the connected characteristic leaves. The
torus preserves every connected refined stratum, so
\eqref{eq:wall-characteristic-relation} is exactly the relation
\eqref{eq:separate-relation} of Definition \ref{def:stratified-system}.
Each leaf is compact, hence closed in $M$.

\emph{The local quotient.}
Let $N$ be an $\mathbf{H}$-invariant normal-form neighborhood of
a point of $M$ with active set $I$. It is saturated for
\eqref{eq:wall-characteristic-relation}, so $p(N)$ is open in $X$ and
$p|_N:N\to p(N)$ is a quotient map. In a normal-form chart, the chosen basis gives
\begin{align*}
 q&=\exp\left(\ii\sum_{j=1}^k\theta_j\mathrm{diag}(b_j)\right),
 \qquad(\theta_1,\ldots,\theta_k)\in(\mathbb R/2\pi\mathbb Z)^k.
\end{align*}
The dual moment coordinates and canonical forms are
\begin{align*}
 \rho_j(\lambda)&=\langle\lambda,b_j\rangle
   =\sum_{i=1}^n\lambda_i b_{j,i},
 &
 \Theta_{\mathbf{H}}&=\sum_{j=1}^k\rho_j\,d\theta_j,
 &
 \omega_{T^*\mathbf{H}}&=\sum_{j=1}^k d\rho_j\wedge d\theta_j.
\end{align*}
These signs agree with \eqref{eq:canonical}.
For $j\le r$, translation by $c_j$ gives
$d\rho_j\wedge d\theta_j=dt_j\wedge d\theta_j$.

For $1\le j\le r$ replace $(t_j,\theta_j)$ by
\begin{align*}
 (x_j,y_j,t_j)&=(t_j\cos\theta_j,t_j\sin\theta_j,t_j).
\end{align*}
Keep the other torus, moment and slice coordinates. The fibers of
this map are exactly the classes of
\eqref{eq:wall-characteristic-relation}, so it identifies the local
quotient topology with the product of cones described above. The
factors with $t_j\ne0$ are smooth cylinders, which gives
\eqref{eq:wall-cone-local-model} at a point with active set $I$.

\emph{Dimensions and ranks.}
Let $S$ be a refined source stratum with active set $I$, lying in
the original stratum $\mathcal{M}_\alpha$ of dimension $d$ and over the
sign piece $B_\beta$. The piece $B_\beta$ is open in the affine subspace
$\{t_j=0:j\in I\}$ of the trace-zero hyperplane, which has codimension
$c$ because $b_1,\ldots,b_k$ is a basis. By \eqref{eq:pointwisespec},
$h$ is a submersion on $\mathcal{M}_\alpha$, so $\dim S=d-c$. The
leaves in $S$ are $K_I$-orbits of dimension $c$, so $\dim p(S)=d-2c$.
Put $\mathcal{P}_\delta=\psi(\mathcal{M}_\alpha)$ and let $T$ be the
refined middle stratum with $S\subseteq\psi^{-1}(T)$. By
Proposition \ref{prop:pointwise-free},
$\psi^{-1}(\mathcal{P}_\delta)=\mathcal{M}_\alpha$ and
$\psi:\psi^{-1}(T)\to T$ is a locally trivial principal
$\mathbf{H}$-bundle with connected base and fiber. Hence $\psi^{-1}(T)$
is connected and equals $S$, and $T=\psi(S)=S/\mathbf{H}$ has dimension
$d-c-k$. By the same proposition, $\pi|_{\mathcal{P}_\delta}$ has
constant rank $k$ and the connected components of its fibers are the
symplectic leaves of $\mathcal{P}_\delta$. Hence $\pi|_T$ is a
submersion into $B_\beta$, of rank $k-c$. Since $T$ is a union of
connected components of fibers of $\pi|_{\mathcal{P}_\delta}$, the
bivector $\Pi_{\mathcal{P}_\delta}$ is tangent to $T$, the bivector
$\Pi_T=\Pi_{\mathcal{P}_\delta}|_T$ has rank $d-2k$, and the connected
components of the fibers of $\pi|_T$ are the symplectic leaves of $T$.
This proves \eqref{eq:wall-rank-formulas}. The fibers of $p$ on $S$ are
$c$-dimensional tori, so $p$ is not injective on $S$ when $c\ge1$, and
the rank of the chain $p(S)\to T\to B_\beta$ is
$\mathrm{rank}\,d(\pi|_T)=k-c$.
\end{proof}

\begin{proof}[Proof of Theorem \ref{thm:integral-wall-quotient}]
\emph{Step 1: the refinements and the relation.}
By Proposition \ref{prop:wall-local-model}, the refined source and
middle strata are finite Whitney refinements satisfying the frontier
condition, the middle refinement is Poisson, each
$\tau_S=\omega|_S$ is presymplectic, and
\eqref{eq:wall-characteristic-relation} is the characteristic relation
\eqref{eq:separate-relation} of Definition \ref{def:stratified-system}.
Every refined stratum is $\mathbf{H}$-invariant, since it is a
connected component of an $\mathbf{H}$-invariant set and $\mathbf{H}$
is connected.

\emph{Step 2: the quotient topology.}

Write $I_\beta=I(\lambda)$ for $\lambda\in B_\beta$. The sign pieces satisfy
\begin{equation}\label{eq:wallsubinci}
 B_{\beta'}\subseteq\overline{B_\beta}
 \quad\Longrightarrow\quad I_\beta\subseteq I_{\beta'}
 \quad\Longrightarrow\quad K_{I_\beta}\subseteq K_{I_{\beta'}}.
\end{equation}
For a closed subset $F\subseteq M$, its saturation is
\[
 p^{-1}(p(F))=\bigcup_\beta
 \bigl(h^{-1}(B_\beta)\cap K_{I_\beta}\cdot F\bigr).
\]
Each $K_{I_\beta}\cdot F$ is closed by compactness. A limit of points
in the $\beta$ summand lies over some
$B_{\beta'}\subseteq\overline{B_\beta}$ and belongs to the $\beta'$ summand
by \eqref{eq:wallsubinci}. Since there are finitely many sign pieces,
the saturation is closed, so $p$ is closed. Its fibers are compact torus
orbits; the standard closed-quotient argument gives Hausdorffness of $X$
and properness of $p$ \cite[Lemma 2.3]{GJSImplosion}.\footnote{For Godement's
smooth quotient criterion, see P.~Petersen,
\emph{Manifolds, Transversality, and de Rham Cohomology}, Theorem~3.3.15,
\url{https://www.math.ucla.edu/~petersen/manifolds.pdf}.}
In particular, the characteristic relation is closed.
The sets $X\setminus p(M\setminus U)$, for finite unions $U$ from a
countable relatively compact basis of $M$, form a relatively compact
basis of $X$. Thus $X$ is second countable, locally compact and paracompact.

\emph{Step 3: smooth charts and the Whitney stratification.}
The local quotients of Proposition \ref{prop:wall-local-model}
cover $X$. Two original normal-form charts are $\mathbf{H}$-equivariant and
preserve $h$. Their overlap map consequently has the form
\begin{align*}
 (\theta,\lambda,y)\longmapsto
 (\theta+f(\lambda,y),\lambda,F(\lambda,y)).
\end{align*}
where $f$ is torus-valued. With $w_j=x_j+\ii y_j$, the induced
transition is
$w'_j=e^{\ii f_j(\lambda,y)}w_j$ and
$y'=F(\lambda,y)$, with the usual translations on the remaining
torus factors. Smoothness of the original quotient charts gives local ambient
smooth extensions of these coefficients and of $F$, and the same
holds for the inverse transition. So these are compatible singular
charts. They define $C_X^\infty$, and $p$ is smooth. On each refined stratum
$S$, the map $p|_S$ forgets the angles $\theta_j$, $j\in I$, so it is a
submersion onto the manifold $p(S)$. 
The images $p(S)$ of connected refined source strata are connected
and pairwise disjoint. The product models identify them with
product strata of the cones, cylinders and $Y$, giving the
Whitney conditions.

For the global frontier condition, every $S$ is
$\mathbf{H}$-invariant. Its closure in $M$ is therefore
$\mathbf{H}$-invariant and saturated for
\eqref{eq:wall-characteristic-relation}.
Since $p$ is continuous and closed,
\begin{align*}
 \overline{p(S)}^{\,X}&=p(\overline S^{\,M}), \qquad  p^{-1}(\overline{p(S)}^{\,X}) =\overline S^{\,M}.
\end{align*}
If $p(S')$ meets $\overline{p(S)}^{\,X}$, then $S'$ meets
$\overline S^{\,M}$. The frontier condition for the refined strata of $M$ gives
$S'\subseteq\overline S^{\,M}$ and hence
$p(S')\subseteq\overline{p(S)}^{\,X}$.
Thus the strata of $X$ satisfy the frontier condition in the
quotient topology.

\emph{Step 4: reduced forms and Poisson brackets.}
The closed form $\omega|_S$ has kernel tangent to $K_I$ and descends to the unique symplectic form on $p(S)=S/K_I$ satisfying $(p|_S)^*\underline\omega_{p(S)}=\omega|_S$. In the local product it is the sum of $dx_j\wedge dy_j/t_j$ on the open cone sheets, the remaining cotangent forms, and the reduced stratal form on $Y$.

Take the product of the linear cone bivectors \eqref{eq:cone-symplectic-form}, the canonical bivectors, and the reduced bivector on $Y$. Invariant polynomial generators for the slice give smooth ambient coefficients for the last factor, as in Lemma \ref{lem:poisson-ranks}. This product bivector is tangent to each stratum and inverts its reduced form. Uniqueness of the reduced forms makes the chart brackets agree, hence they define a stratified Poisson structure on $C_X^\infty$. The local formulas also make $p$ Poisson from the original bracket on $M$. Equivalently, this identity holds off the walls and extends by continuity. Thus the reduced forms are Poisson compatible.

\emph{Step 5: the factored chain.}
The map $\psi$ is constant on every $K_I$-orbit. In the charts
above, the induced map $\underline\psi$ forgets the angles and the
coordinates $x_j,y_j$, and keeps all spectral coordinates $\lambda$,
including the $t_j$, and the slice point $y$. So $\underline\psi$ is a
smooth, stratified Poisson map to $P$, and the diagram
\eqref{eq:global-reduction-diagram} commutes. Let $S$ be a refined
source stratum with active set $I$, over the refined middle stratum
$T=\psi(S)$ and the sign piece $B_\beta$. The torus $\mathbf{H}/K_I$
acts freely and in a Hamiltonian way on $p(S)=S/K_I$, and
$\underline\psi|_{p(S)}:p(S)\to T$ is its quotient map, a surjective
submersion with connected fibers. By
Proposition \ref{prop:wall-local-model},
$\dim T=\dim p(S)-(k-c)$, the map $\pi|_T$ has constant rank $k-c$,
and the connected components of its fibers are the symplectic leaves
of $T$. Hence the bottom row of \eqref{eq:stratum-reduction-diagram}
is a strong chain of rank $k-c$. Finally, for compact $D\subset P$ the
set $\underline\psi^{-1}(D)=p(\psi^{-1}(D))$ is compact, since $\psi$
is proper. This proves stratified superintegrability in the sense of
Definition \ref{def:stratified-system}. The spectral torus action
descends to a smooth action on $X$ which rotates the cone
coordinates, so the complete Hamiltonian flows of
Proposition \ref{prop:pointwise-free} descend to complete flows
preserving every stratum of $X$. On a refined stratum $S$ lying on $c$
walls, the flow of a spectral Hamiltonian $F(h)$ descends to a flow of
the residual torus $\mathbf{H}/K_I$ on $p(S)$, since the $c$ directions
of $K_I$ are collapsed by $p$. When $c=k$ the descended spectral
Hamiltonians are constant, which is the rank-zero case allowed in
Definition \ref{def:strong}.
\end{proof}

\begin{remark}[Scope of the wall theorem]
\label{rem:integral-wall-scope}
An element $b\in\Lambda$ is \emph{primitive} if $b\ne0$ and $b=q b'$ with $b'\in\Lambda$, $q\in\mathbb Z$, implies $|q|=1$.
The theorem covers any single wall with a primitive integral normal,
and any family of walls whose normals can be extended to a basis of
the lattice. For example, any set of walls
$\lambda_i-\lambda_{i+1}=a_i$ given by simple roots. Every wall
$\lambda_I-\lambda_J=\mu_K$ coming from the complementary Horn
condition also has a primitive normal, since its nonzero coordinates
are $1$ and $-1$. The quotient contains the free parts of the chosen
walls together with the regions on both sides. Nonfree points and
repeated spectra stay outside the domain.

The cone atlas determines the smooth functions on $X$. These are not
all continuous functions whose pullback along $p$ is smooth. For instance,
$\overline f=(x^2-y^2)/t$ on the open sheets, extended by $0$ at the
vertex, is continuous with smooth pullback $s\cos(2\theta)$, but it
has no ambient extension of class $C^1$ at the vertex: the derivative
of such an extension along the rays $\rho(\cos\theta,\sin\theta,1)$,
$\rho>0$, would be $A\cos\theta+B\sin\theta+D$ for constants
$A,B,D$, whereas it equals $\cos(2\theta)$.
\end{remark}

\section{\texorpdfstring{The $\mathrm{SU}(3)$ system and topological obstructions}{The SU(3) system and topological obstructions}}
\label{sec:su3}\label{sec:examples}\label{sec:su3-calculation}

For the orbit with spectrum $(a,0,-a)$, condition
\eqref{hyp:spectrum} excludes two interior walls, although the free
set still has points over these walls and over their intersection. We
compute both domains and apply Theorem \ref{thm:integral-wall-quotient}.
We then show, for the zero orbit, that the leaf relation cannot be
extended to the whole phase space.

\subsection{The allowed spectral domain for \texorpdfstring{$\mathrm{SU}(3)$}{SU(3)}}\label{subsec:su3-condition}

 Fix $\xi=\ii\mathrm{diag}(a,0,-a)$, where $a>0$. Write
\begin{align}\label{eq:su3-gap-coordinates} u=\lambda_1-\lambda_2,\qquad v=\lambda_2-\lambda_3,\qquad \lambda=\frac13(2u+v,v-u,-u-2v). \end{align}
These are the coordinates $u,v$ of \cite[Sections 5.3.1 and 5.5]{ChenEtAl}. In particular, $x$ is conjugate to $\ii(u\vec\omega_1+v\vec\omega_2)$, where $\vec\omega_1=\frac13\mathrm{diag}(2,-1,-1)$ and $\vec\omega_2=\frac13\mathrm{diag}(1,1,-2)$.

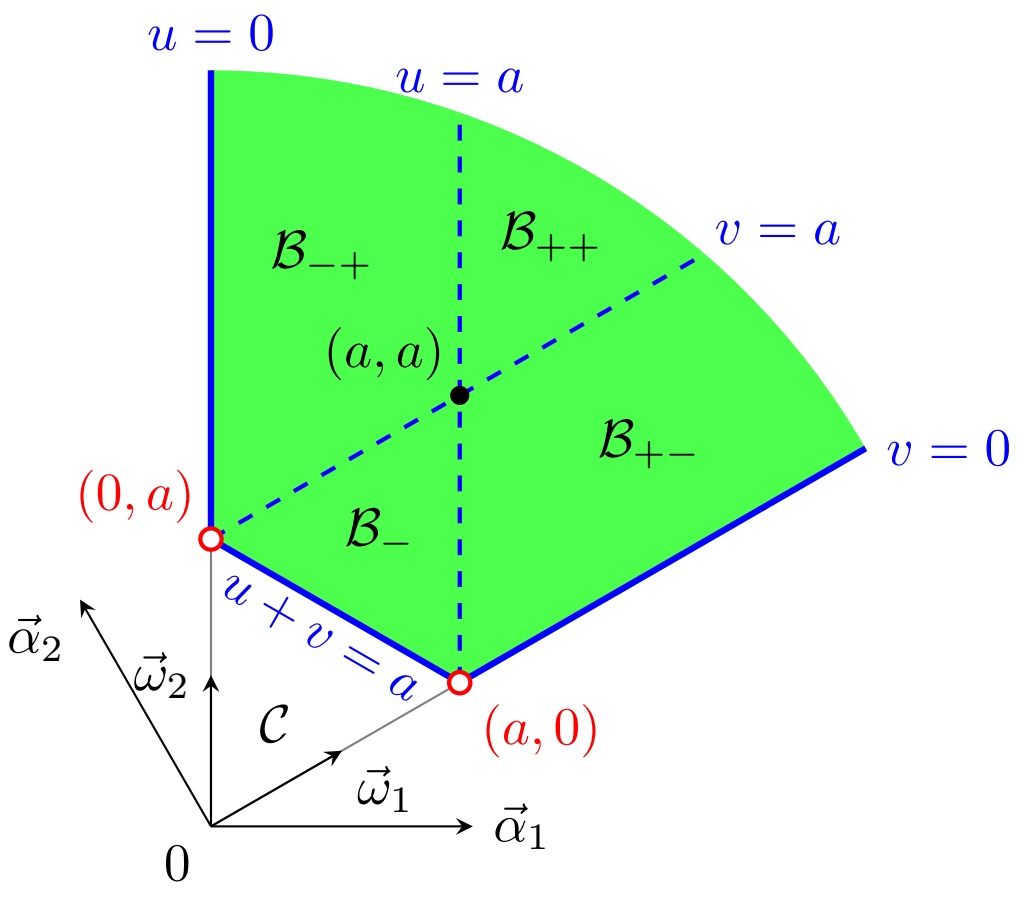
\begin{figure}[H]
\centering
\begin{tikzpicture}[>=stealth,font=\small]
\def\orbA{1.9}  
\def\cutR{5.0}  
\pgfmathsetmacro{\wallend}{(-\orbA+sqrt(4*\cutR*\cutR-3*\orbA*\orbA))/2} 
\begin{scope}[x={(0.8660cm,0.5cm)},y={(0cm,1cm)}]
  \coordinate (O)   at (0,0);
  \coordinate (w1)  at (1,0);
  \coordinate (w2)  at (0,1);
  \coordinate (al1) at (2,-1);
  \coordinate (al2) at (-1,2);
  \coordinate (A1)  at (\orbA,0);
  \coordinate (A2)  at (0,\orbA);
  \coordinate (AA)  at (\orbA,\orbA);
  \coordinate (M)   at ({\orbA/2},{\orbA/2});
  \coordinate (E1)  at (\orbA,\wallend);
  \coordinate (E2)  at (\wallend,\orbA);
  \coordinate (Lm)  at ({0.68*\orbA},{0.68*\orbA});
  \coordinate (Lmp) at (0.85,3.35);
  \coordinate (Lpm) at (3.35,0.85);
  \coordinate (Lpp) at (2.6,2.6);
\end{scope}
\fill[green!70] (A2) -- (A1) -- (30:\cutR) arc[start angle=30,end angle=90,radius=\cutR] -- cycle;
\begin{scope}
  \clip (A2) -- (A1) -- (30:\cutR) arc[start angle=30,end angle=90,radius=\cutR] -- cycle;
  \begin{scope}[x={(0.8660cm,0.5cm)},y={(0cm,1cm)}]
    \draw[blue,dashed,thick] (\orbA,0) -- (\orbA,6);
    \draw[blue,dashed,thick] (0,\orbA) -- (6,\orbA);
  \end{scope}
\end{scope}
\draw[blue,very thick] (90:\cutR) -- (A2) -- (A1) -- (30:\cutR);
\draw[gray] (O) -- (A1);
\draw[gray] (O) -- (A2);
\draw[->] (O) -- (al1) node[right] {$\vec\alpha_1$};
\draw[->] (O) -- (al2) node[below left=-1pt] {$\vec\alpha_2$};
\draw[->] (O) -- (w1)  node[below right=-1pt] {$\vec\omega_1$};
\draw[->] (O) -- (w2)  node[left]  {$\vec\omega_2$};
\node[below left] at (O) {$0$};
\node at (58:0.8) {$\mathcal C$};
\draw[red,thick,fill=white] (A1) circle (2pt) node[below right,text=red] {$(a,0)$};
\draw[red,thick,fill=white] (A2) circle (2pt) node[above left=-1pt,text=red] {$(0,a)$};
\fill (AA) circle (1.8pt) node[above left=-1pt] {$(a,a)$};
\node[blue,above] at (E1) {$u=a$};
\node[blue,above right=-2pt] at (E2) {$v=a$};
\node[blue,right] at (30:\cutR) {$v=0$};
\node[blue,above] at (90:\cutR) {$u=0$};
\node[blue,rotate=-30,anchor=north,inner sep=2pt] at (M) {$u+v=a$};
\node at (Lm)  {$\mathcal B_{-}$};
\node at (Lmp) {$\mathcal B_{-+}$};
\node at (Lpm) {$\mathcal B_{+-}$};
\node at (Lpp) {$\mathcal B_{++}$};
\end{tikzpicture}
\caption{The feasible set $\mathcal{B}_\xi$ of \eqref{eq:su3-feasible} in the
closed chamber $\mathcal{C}$, with $(u,v)$ drawn at
$u\vec\omega_1+v\vec\omega_2$; here $\vec\alpha_1=2\vec\omega_1-\vec\omega_2$
and $\vec\alpha_2=2\vec\omega_2-\vec\omega_1$. Only a bounded part of
$\mathcal{B}_\xi$ (green) is shown. The curved outer edge only truncates the drawing. Condition \eqref{hyp:spectrum} removes
the actual boundary portions $u=0$, $v=0$, and $u+v=a$ and the
dashed walls $u=a$, $v=a$, leaving the four
components \eqref{eq:su3-four-components} of $\Omega_\xi$. Compare
\cite[Figure 5]{ChenEtAl}, where the walls are the lines
$\mathcal{C}_{(1)}$ and $\mathcal{C}_{(2)}$.}
\label{fig:su3-domain}
\end{figure}

Recall that $\mathcal{B}_\xi=h(\mathcal{R}_\xi)$ is the set of feasible spectra. Every point of $\mathcal{R}_\xi$ has a representative $(x,\gamma)$ with $x=\ii D$ and $D=\mathrm{diag}(\lambda_1,\lambda_2,\lambda_3)$, and then $-\ii J(x,\gamma)=D-\gamma^*D\gamma$. Thus $\lambda\in\mathcal{C}$ is feasible exactly when
\[
D-\gamma^*D\gamma\ \text{has spectrum}\ (a,0,-a)\ \text{for some}\ \gamma\in\mathrm{SU}(3).
\]
In the coordinates \eqref{eq:su3-gap-coordinates}, this set is
\begin{align}\label{eq:su3-feasible} \mathcal{B}_\xi=\{(u,v):u\ge0,\ v\ge0,\ u+v\ge a\}, \end{align}
the region drawn in Figure \ref{fig:su3-domain}; see also \cite[Lemma 5.8]{ChenEtAl}. Indeed, the largest eigenvalue $a$ of $D-\gamma^*D\gamma$ is at most $\lambda_1-\lambda_3=u+v$. Conversely, if $u,v\ge0$ and $d=u+v\ge a$, put $s=a/d$, $c=(1-s^2)^{1/2}$ and
\begin{align}\label{eq:su3-feasibility-matrices} \gamma=\begin{pmatrix}c&0&s\\0&1&0\\-s&0&c\end{pmatrix}\in\mathrm{SU}(3). \end{align}
Then $D-\gamma^*D\gamma$ vanishes on $e_2$ and is trace-free on $\mathrm{span}(e_1,e_3)$ with square $a^2I_2$, so its spectrum is $(a,0,-a)$.

On $\mathcal{B}_\xi$, condition \eqref{hyp:spectrum} is equivalent to
\begin{align}\label{eq:su3-H-explicit} u>0,\qquad v>0,\qquad u+v>a,\qquad u\ne a,\qquad v\ne a. \end{align}
Here $u,v>0$ is the condition $\lambda\in\mathcal{C}^\circ$. In dimension three, every decomposition in \eqref{hyp:spectrum} has blocks of sizes one and two, and taking complements interchanges them. So it suffices to take $I=\{i\}$, $J=\{j\}$ and $K=\{k\}$ with $i\ne j$. The $1\times1$ problem is solvable exactly when $\lambda_i-\lambda_j=\mu_k$. Since $|\lambda_i-\lambda_j|\in\{u,v,u+v\}$ is nonzero, this forces $\mu_k=\pm a$ and one of $u$, $v$, $u+v$ to equal $a$. Conversely, if $u$, $v$ or $u+v$ equals $a$, take $\{i,j\}=\{1,2\}$, $\{2,3\}$ or $\{1,3\}$, respectively. The complementary $2\times2$ blocks have eigenvalue gaps $(v,u+v)$, $(u,u+v)$ or $(u,v)$, and the target gap is $a$. The gap test \eqref{eq:two-dimensional-gap-test} holds in each case, and the trace identity follows from the total trace. Both block problems are then solvable, so \eqref{hyp:spectrum} fails.

Hence $\Omega_\xi$ is the interior of $\mathcal{B}_\xi$ with the two walls $u=a$ and $v=a$ removed. Its connected components are the four convex sets
\begin{align}\label{eq:su3-four-components} \mathcal{B}_{-}&=\{0<u<a,\ 0<v<a,\ u+v>a\},&\mathcal{B}_{-+}&=\{0<u<a,\ v>a\},\\ \mathcal{B}_{+-}&=\{u>a,\ 0<v<a\},&\mathcal{B}_{++}&=\{u>a,\ v>a\},\nonumber \end{align}
since a path joining two of them must cross a wall.

\begin{remark}\label{rem:su3-two-conditions}
For this particular orbit, \eqref{hyp:spectrum} and \eqref{hyp:trace} agree on $\mathcal{B}_\xi$. Indeed, for $I\ne J$ with $|I|=|J|\in\{1,2\}$, the differences $\lambda_I-\lambda_J$ are $\pm u$, $\pm v$, or $\pm(u+v)$. The proper subset sums $\mu_K$ are $0$, $a$, and $-a$. Condition \eqref{hyp:trace} is therefore exactly \eqref{eq:su3-H-explicit}. Proposition \ref{prop:comparison} shows that the implication is strict for the different orbit $\mu=(3,-1,-2)$.
The present example makes the complete-fiber domain and the wall construction explicit; the strictness example concerns a different orbit.
\end{remark}

\subsection{Free points in excluded complete fibers}\label{subsec:su3-restriction-effect}

For a feasible simple spectrum $\lambda\notin\Omega_\xi$, the fiber
$h^{-1}(\lambda)$ contains a point at which the spectral torus does
not act freely. The same fiber may also contain free points.

For example, at $(u,v)=(a,2a)$ put $D=(a/3)\mathrm{diag}(4,1,-5)$ and $x=\ii D$. The rotation \eqref{eq:su3-feasibility-matrices}, with $s=1/3$, gives
\begin{align*}
 G_{x,\gamma}=G_{x,y}=\{\mathrm{diag}(t,t^{-2},t):|t|=1\},
\end{align*}
so its reduced point is free for the spectral torus. In contrast,
\begin{align*}
 \widehat\gamma=\begin{pmatrix}0&1&0\\-1&0&0\\0&0&1\end{pmatrix},
 \qquad D-\widehat\gamma^*D\widehat\gamma=\mathrm{diag}(a,-a,0)
\end{align*}
has $G_{x,y}=\mathbf{H}$ and $G_{x,\widehat\gamma}=\{\mathrm{diag}(t,t,t^{-2}):|t|=1\}$, hence a one-dimensional spectral-torus stabilizer. The two representatives have the same conjugation orbit type, since
their stabilizer circles are conjugate, but different spectral-torus
stabilizers. The restriction to $\Omega_\xi$ removes both
points because it retains only complete free fibers.

\subsection{Retained strata and how they meet}\label{subsec:su3-retained-geometry}
The space $\mathcal{M}_\xi=h^{-1}(\Omega_\xi)$ has symplectic strata of dimensions six and four. We classify them and exhibit an incidence at fixed spectrum.

\subsubsection{The two orbit types and their Poisson chains}\label{subsec:su3-retained-strata}

The two stabilizer types are:
\begin{align*}
Z_3=\{tI_3:t^3=1\},\qquad H(2,1) =\{\mathrm{diag}(t,t,t^{-2}):|t|=1\}.
\end{align*}
We use the notation of \cite[Notation 5.2]{ChenEtAl}. The dimensions and ranks are summarized below. All dimensions are real. Here $\mathcal{M}_\alpha$ is a connected orbit-type stratum of $\mathcal{M}_\xi$ and $\mathcal{P}_\delta=\psi(\mathcal{M}_\alpha)$ is its image.
\begin{center}
\renewcommand{\arraystretch}{1.2}
\begin{tabular}{c|ccccc}
Stabilizer type & $\dim\mathcal{M}_\alpha$ & $\dim\mathcal{P}_\delta$ & $\mathrm{rank}\Pi_\delta$ & $\mathrm{rank}d\pi$ & fiber of $\psi$\\\hline
$Z_3$ & $6$ & $4$ & $2$ & $2$ & $\mathbb T^2$\\
$H(2,1)$ & $4$ & $2$ & $0$ & $2$ & $\mathbb T^2$
\end{tabular}
\end{center}
Here $\mathbb T^2=(S^1)^2$.

Diagonalize $x=\ii D$, where $D=\mathrm{diag}(\lambda_1,\lambda_2,\lambda_3)$ and $\lambda_1>\lambda_2>\lambda_3$. Every element of $G_x$ is diagonal, say $g=\mathrm{diag}(t_1,t_2,t_3)$. To commute also with $\gamma$, it must satisfy
\begin{align*}
(t_i-t_j)\gamma_{ij}=0\qquad(1\le i,j\le3).
\end{align*}
These equations can be read from a graph with vertices $1,2,3$: join $i$ and $j$ when $\gamma_{ij}\ne0$ or $\gamma_{ji}\ne0$. The numbers $t_i$ must be equal on each connected component.

A connected graph forces $t_1=t_2=t_3$, giving $G_{x,\gamma}=Z_3$. A graph with one two-vertex component and one isolated vertex gives a stabilizer conjugate to $H(2,1)$. Three isolated vertices would make $\gamma$ diagonal, hence $J(x,\gamma)=0$, which is impossible for $a>0$. This proves the classification. For the stabilizers $G_{x,y}$ of the pairs $(x,y)$, Proposition \ref{prop:stabilizer} gives
\begin{align*}
G_{x,y}=G_{x,\gamma}
\end{align*}
throughout $\mathcal{M}_\xi$, so the same two stabilizer types describe the classes $[(x,y)]\in\mathcal{P}_\xi^\Omega$.

\noindent\textbf{Where the four-dimensional strata occur.}
For type $H(2,1)$, let $i$ be the isolated index and let $p<q$ be the other two indices. Then $\gamma$ preserves the line $\C e_i$ and has a $2\times2$ unitary block $U$ on $\mathrm{span}(e_p,e_q)$. Put
\begin{align*}
d=\lambda_p-\lambda_q,\qquad s=|\gamma_{pq}|.
\end{align*}
The isolated line contributes the zero eigenvalue of $D-\gamma^*D\gamma$. The other two eigenvalues are $\pm ds$, so the fixed-orbit condition is simply $ds=a$. Since $s\le1$ and the equality $d=a$ is excluded by \eqref{eq:su3-H-explicit}, the allowed gaps are
\begin{align}\label{eq:su3-K-domains}
ds=a,\qquad d>a,\qquad
(i,d)=(1,v),\ (2,u+v),\ \text{or }(3,u).
\end{align}
Conversely, whenever $d>a$, take $s=a/d$ and $c=\sqrt{1-s^2}$. The block
\begin{align*}
U=\begin{pmatrix}c&s\\-s&c\end{pmatrix}
\end{align*}
and the identity on $\C e_i$ give a representative. In particular, each entry of this block is nonzero.

At fixed $i$ and $\lambda$, there is exactly one class $[(x,y)]\in\mathcal{P}_\xi^\Omega$ of this type. Indeed, the diagonal entries of $U^*\mathrm{diag}(\lambda_p,\lambda_q)U$ are $\lambda_p-ds^2$ and $\lambda_q+ds^2$, while the modulus of its off-diagonal entry is $dsc$. These quantities depend only on $\lambda$, and diagonal conjugation removes the remaining phase. The displayed rotation varies smoothly with $\lambda$, so the spectral map identifies this set of classes $[(x,y)]$ smoothly with its allowed spectral domain.

This local spectral parametrization gives middle strata
of dimension two. The local product with the two-dimensional torus
$\mathbf H$ in Theorem \ref{thm:main} then gives source strata of
dimension four.

For type $Z_3$, the effective conjugation action is free, so $\dim\mathcal{M}_\alpha=\dim T^*G-2\dim G+\dim\mathcal{O}_\xi=16-16+6=6$. The other type has dimension four by the preceding spectral parametrization and torus bundle. Formula \eqref{eq:main-dimensions}, with $k=2$, gives the remaining table entries. The middle symplectic leaves have dimension two in the first row and are points in the second. Both chains are strong, with characteristic projection $p=\mathrm{id}$.

\subsubsection{How the six- and four-dimensional strata meet}\label{subsubsec:su3-incidence}
The two orbit types are not separated by the spectral map: the
following path keeps $x$, and hence the spectrum, fixed, and moves from
a six-dimensional stratum into a four-dimensional one by varying
$\gamma$ alone.

\begin{proposition}  \label{prop:singular-example}
Let $a>0$, $G=\mathrm{SU}(3)$, and $\xi=\ii\mathrm{diag}(a,0,-a)$.
Then there exist connected symplectic orbit-type strata $S_6,S_4\subset\mathcal{M}_\xi$ of dimensions six and four, respectively, a number $r>0$, and a continuous path $m:[0,r)\to\mathcal{M}_\xi$ such that
\begin{align*}
 h(m(t))&=(3a,0,-3a) \text{ with } m(0) \in S_4,\qquad m(t)\in S_6\quad 0<t<r .
\end{align*}
\end{proposition}
\begin{proof}

Put $L=3a$, $c=L-a^2/(2L)$, $r=(L^2-c^2)^{1/2}$ and $s(t)=(r^2-t^2)^{1/2}$. Set
\begin{equation}\label{eq:singular-path}
A=\begin{pmatrix}0&L&0\\L&0&0\\0&0&0\end{pmatrix},\quad
\gamma(t)=\begin{pmatrix}1&0&0\\0&(c-\ii s(t))/L&t/L\\0&-t/L&(c+\ii s(t))/L\end{pmatrix},\quad m(t)=[(\ii A,\gamma(t))].
\end{equation}
The lower block has determinant one and is unitary because $c^2+s(t)^2+t^2=L^2$. Moreover
\[
\gamma(t)^*A\gamma(t)=\begin{pmatrix}0&c-\ii s(t)&t\\c+\ii s(t)&0&0\\t&0&0\end{pmatrix}.
\]
Its difference from $A$ has eigenvalues $0$ and $\pm\sqrt{(L-c)^2+s(t)^2+t^2}=\pm a$. Thus $J\in\mathcal{O}_\xi$, while the constant spectrum $(L,0,-L)$ lies in $\Omega_\xi$.

At $t=0$, $\gamma(0)=\mathrm{diag}(1,e^{-\ii\theta},e^{\ii\theta})$ with $\cos\theta=17/18$. Its eigenvalues are distinct. A common centralizer with $A$ must be diagonal with equal first two entries, giving $H(2,1)$. For $t>0$, every matrix commuting with $A$ has form
$g=\left(\begin{smallmatrix}\alpha&\beta&0\\\beta&\alpha&0\\0&0&\delta\end{smallmatrix}\right)$.
The $(1,3)$ and $(2,3)$ entries of $g\gamma(t)=\gamma(t)g$ imply $\beta t=0$ and $(\alpha-\delta)t=0$, hence $g\in Z_3$. The dimension table and connectedness of $(0,r)$ now give $S_4,S_6$ and the stated limit.

\end{proof}

The change occurs in the stabilizer, while both the coadjoint orbit and the conserved spectrum stay fixed. In particular, a four-dimensional stratum meets the closure of a six-dimensional one inside a single spectral fiber. 

\subsection{The two-wall quotient}
We now apply Theorem \ref{thm:integral-wall-quotient} to the two
interior walls $u=a$ and $v=a$ excluded by condition
\eqref{hyp:spectrum}, on the whole free domain of the nonzero $\mathrm{SU}(3)$
orbit. The first step is to identify that domain. The second is to
read off the subtori and the dimension table from the theorem.

\begin{corollary}[The two interior walls for $\mathrm{SU}(3)$]
\label{cor:su3-wall-quotient}
Let $\xi=\ii\mathrm{diag}(a,0,-a)$, where $a>0$, and use the
spectral coordinates $u=\lambda_1-\lambda_2$ and
$v=\lambda_2-\lambda_3$. Then
\begin{equation}\label{eq:su3-free-base}
 \mathcal{B}_\xi^{\mathrm{free}}
 =D_a:=\{(u,v):u>0,\ v>0,\ u+v>a\}.
\end{equation}
The complete-fiber domain is
\[
\Omega_\xi=D_a\setminus\bigl(\{u=a\}\cup\{v=a\}\bigr).
\]
The wall quotient therefore retains the free points on the two
interior walls omitted by the complete-fiber restriction.
On the whole space $\mathcal{M}_\xi^{\mathrm{free}}$, simultaneously
refine by the natural walls $u=a$ and $v=a$. The base has nine sign
pieces: four open chambers, the four connected wall components obtained by
removing $(a,a)$ from the walls, and the point $(a,a)$.

The construction of Theorem \ref{thm:integral-wall-quotient} gives one
proper Hausdorff characteristic quotient
\[
 p:\mathcal{M}_\xi^{\mathrm{free}}\longrightarrow X_a.
\]
Its fibers are points off the walls, circles over each wall component,
and two-dimensional tori over $(a,a)$. These are exactly the connected
characteristic leaves of the refined source strata. The quotient has
the Whitney stratification, cone smooth structure, and
Poisson structure of that theorem, and the factorization
\[
 \mathcal{M}_\xi^{\mathrm{free}}
 \xrightarrow{\ p\ }X_a
 \xrightarrow{\ \underline\psi\ }\mathcal{P}_\xi^{\mathrm{free}}
 \xrightarrow{\ \pi\ }D_a
\]
is the global characteristic quotient of the chain
$\mathcal{M}_\xi^{\mathrm{free}}\xrightarrow{\psi}\mathcal{P}_\xi^{\mathrm{free}}\xrightarrow{\pi}D_a$
with these refinements. All original orbit-type strata of
$\mathcal{M}_\xi^{\mathrm{free}}$ are kept, both the six-dimensional
ones of type $Z_3$ and the four-dimensional ones of type $H(2,1)$.

For every nonempty refined stratum, let $d$ be the dimension of the
original symplectic stratum containing it and $c$ the number of walls
containing its image in the base. The dimensions and strong-chain ranks are as follows:
\begin{center}
\renewcommand{\arraystretch}{1.15}
\begin{tabular}{cc|cccc}
$d$ & $c$ & $\dim S=d-c$ & $\dim p(S)=d-2c$ & $\dim\psi(S)=d-2-c$ & rank $2-c$\\\hline
$6$ & $0$ & $6$ & $6$ & $4$ & $2$\\
$6$ & $1$ & $5$ & $4$ & $3$ & $1$\\
$6$ & $2$ & $4$ & $2$ & $2$ & $0$\\
$4$ & $0$ & $4$ & $4$ & $2$ & $2$\\
$4$ & $1$ & $3$ & $2$ & $1$ & $1$\\
$4$ & $2$ & $2$ & $0$ & $0$ & $0$
\end{tabular}
\end{center}
\end{corollary}

\begin{proof}
By \eqref{eq:su3-feasible}, the feasible spectra with distinct
entries form the set $\{u>0,v>0,u+v\geq a\}$. On the free set, $h$ is
a submersion on each original symplectic stratum, so its image is open
in the open chamber and cannot contain the boundary line $u+v=a$.

We construct, over every $(u,v)\in D_a$, a free point of orbit type $Z_3$.
Let $D=\mathrm{diag}\lambda$ and choose the unique
$\tau\in(0,1/2)$ satisfying
\begin{equation}\label{eq:su3-reflection-parameter}
 (u^2+v^2)\tau-(u-v)^2\tau^2=\frac{a^2}{4}.
\end{equation}
The left side increases strictly on $[0,1/2]$ from zero to
$(u+v)^2/4>a^2/4$. Put
\[
 e=\bigl(\sqrt\tau,\sqrt{1-2\tau},\sqrt\tau\bigr)^T,
 \qquad \gamma=2ee^T-I_3\in\mathrm{SO}(3)\subset\mathrm{SU}(3),
 \qquad x=\ii D.
\]
Writing $b=e^TDe$ and $w=(D-bI_3)e$, we have $e^Tw=0$ and
\[
 D-\gamma D\gamma=2(we^T+ew^T),\qquad
 \|w\|^2=(u^2+v^2)\tau-(u-v)^2\tau^2=\frac{a^2}{4}.
\]
Thus $D-\gamma D\gamma$ has eigenvalues $a,0,-a$, as required.
All off-diagonal entries of $\gamma$ are nonzero. For $i\ne j$,
\[
 (\gamma D\gamma)_{ij}
 =2e_ie_j(2b-\lambda_i-\lambda_j).
\]
The three pair sums of the distinct $\lambda_i$ are distinct, so at
most one of these three off-diagonal entries vanishes. The graphs of nonzero off-diagonal entries for $\gamma$ and
$\gamma D\gamma$, as defined in Section \ref{subsec:su3-retained-strata}, are therefore
connected. Consequently
\[
 G_{x,\gamma}=G_{x,-\gamma^{-1}x\gamma}=Z_3.
\]
This proves freeness and the reverse inclusion in
\eqref{eq:su3-free-base}.

The representative $\gamma=2ee^T-I_3$ has eigenvalues $1,-1,-1$, so
$\gamma\notin G_{\mathrm{reg}}$, although $x$ has simple spectrum and
the reduced point has orbit type $Z_3$ and is free for the spectral
torus. Corollary \ref{cor:su3-wall-quotient} imposes no eigenvalue
multiplicity condition on $\gamma$.

Since $(a,a)\in D_a$, the preceding construction also
gives a free point of type $Z_3$ at the wall intersection. Thus the
two-dimensional characteristic leaves there occur on the
six-dimensional strata.

The integral normals
$b_1=(1,-1,0)$ and $b_2=(0,1,-1)$ form a basis of the lattice
$\{b\in\mathbb Z^3:\sum_i b_i=0\}$ of $\mathbf{H}$.
Apply Theorem \ref{thm:integral-wall-quotient} with
$t_1=u-a$ and $t_2=v-a$. Its subtori are the circle generated by
$b_1$ on $u=a$, the circle generated by $b_2$ on $v=a$, and all of
$\mathbf{H}$ at the intersection. Its dimension formulas give the
table.

The classification in Section \ref{subsec:su3-retained-strata}
used only that $x$ has simple spectrum and $J\ne0$, so the orbit types
on the free set are still exactly $Z_3$ and $H(2,1)$, and at free
points $G_{x,y}=G_{x,\gamma}$. The four-dimensional free strata have isolated index
$i=1,2,3$ with, respectively, $v>a$, $u+v>a$, and $u>a$: the block
equation is $ds=a$ as before. If $d=a$ the block is off-diagonal and
$y$ is diagonal, so the point is not free, while if $d>a$ the two
stabilizers agree. For the isolated index $2$, the $(1,3)$-rotation
with $\sin\theta=a/(u+v)$ is free throughout $D_a$. Hence the
quotient $X_a$ also contains strata coming from the lower orbit type,
including a two-dimensional refined source stratum at $(a,a)$, contained in an original four-dimensional stratum and reduced by $p$ to a point.
\end{proof}

Away from the two walls, \eqref{eq:su3-H-explicit} shows that this
domain is exactly $\mathcal{M}_\xi$ over the four components of
$\Omega_\xi$, so the new chain extends the chain of
Theorem \ref{thm:main} on all four chambers, with the same orbit
types and the same closure relations between strata there.

\begin{remark}\label{rem:free-nonproper}
On the free domain of Proposition \ref{prop:pointwise-free},
the spectral projection need not be proper.
For $r>0$ sufficiently small, put
\begin{align*}
 &u_r=\frac{a}{\sqrt{1-r^2}},\qquad v_r=\frac{5a-u_r}{2},\\
 &D_r=\mathrm{diag}\!\left(\frac{2u_r+v_r}{3},\frac{v_r-u_r}{3},-\frac{u_r+2v_r}{3}\right),\qquad
 \gamma_r=\begin{pmatrix}r&\sqrt{1-r^2}&0\\-\sqrt{1-r^2}&r&0\\0&0&1\end{pmatrix}.
\end{align*}
The nonzero eigenvalues of $D_r-\gamma_r^*D_r\gamma_r$ are $\pm u_r\sqrt{1-r^2}=\pm a$. Thus $m_r=[(\ii D_r,\gamma_r)]$ belongs to $\mathcal{M}_\xi$, since $u_r>a$ and $v_r>a$. As $r\to0$, it converges in $\mathcal{R}_\xi$ to the nonfree point at spectrum $(u,v)=(a,2a)$ computed in Section \ref{subsec:su3-restriction-effect}. That section also gives a free point at this same spectrum, so $(a,2a)\in\mathcal{B}_\xi^{\mathrm{free}}$. For any sequence $r_\nu\downarrow0$, the spectra of $m_{r_\nu}$ together with their limit form a compact subset of $\mathcal{B}_\xi^{\mathrm{free}}$. Its inverse image under $h$ in $\mathcal{M}_\xi^{\mathrm{free}}$ is not compact: the displayed sequence has a unique limit in the Hausdorff ambient space, and that limit is excluded. The spaces are metrizable, being semialgebraic subsets of Euclidean spaces in invariant polynomial coordinates. Hence $h|_{\mathcal{M}_\xi^{\mathrm{free}}}$ is not proper. The properness of $\psi$ on this domain is proved independently in Proposition \ref{prop:pointwise-free}. Together with $h=\pi\circ\psi$, it implies that $\pi|_{\mathcal{P}_\xi^{\mathrm{free}}}$ cannot be proper either.
\end{remark}

\subsection{The zero orbit: obstruction and spectral restriction}\label{sec:unconditional}
We next examine the characteristic quotient on the unrestricted zero-orbit chain.
We compute it for $\mathrm{SU}(2)$ and prove it for every $n\ge2$.
The simple-spectrum restriction then gives an explicit strong
Liouville chain.

\subsubsection{The complete calculation for \texorpdfstring{$\mathrm{SU}(2)$}{SU(2)}}\label{subsec:zero-su2}

For $\xi=0$ the commuting matrices $x,\gamma$ are simultaneously diagonalizable. In $\mathrm{SU}(2)$ write
\begin{equation}\label{eq:zero-su2-pair}
x_s=\ii\mathrm{diag}(s,-s),\quad
\gamma_\theta=\mathrm{diag}(e^{\ii\theta},e^{-\ii\theta}),\quad s\in\mathbb R,\quad\theta\in\mathbb R/2\pi\mathbb Z.
\end{equation}
The only residual equivalence is $(s,\theta)\sim(-s,-\theta)$, including $s=0$, where the unordered eigenvalues of $\gamma$ determine its class. Consequently
\begin{align}
\mathcal R_0&\cong(\mathbb R\times S^1)/\{\pm1\},\label{eq:zero-su2-quotient}\\
\mathcal R_0&\xrightarrow{\psi}[0,\infty)\xrightarrow{\mathrm{id}}[0,\infty),\qquad \psi([(s,\theta)])=u=2s^2.\label{eq:zero-su2-chain}
\end{align}
The first identification is a homeomorphism: the continuous bijection is proper because $u$ bounds $s$ and $S^1$ is compact. The target coordinate $u=\mathrm{tr}((-\ii x)^2)$ gives its quotient smooth structure. Compare \cite[Section 4.7(I)]{ChenEtAl}.

There is one two-dimensional orbit-type stratum, with stabilizer $\mathbf{H}$, and two point strata $m_+=[(0,I_2)]$, $m_-=[(0,-I_2)]$ with stabilizer $G$. Off these points local diagonalization of either $x$ or
$\gamma$ supplies smooth charts. Refining the targets into $\{0\}$ and $(0,\infty)$ gives the fiber
\begin{equation}\label{eq:zero-su2-interval}
F=\psi^{-1}(0)=\{[(0,\gamma_\theta)]:0\le\theta\le\pi\}.
\end{equation}
Its open interval $F^\circ$ is one refined source stratum with zero restricted form, hence one characteristic leaf. Each endpoint is a different leaf. Since $F^\circ$ is not closed in $\mathcal R_0$, its image under \eqref{eq:separate-relation} is not closed. The quotient is not $T_1$, and no Hausdorff topology on these classes makes the projection continuous.

\subsubsection{The obstruction for every \texorpdfstring{$n\ge2$}{n greater than or equal to 2}}\label{subsec:zero-sun}

The next lemma uses the fact that \eqref{eq:separate-relation} never identifies different refined strata. It gives the obstruction in Proposition \ref{prop:unrestricted-zero} and Section \ref{app:disk}.

\begin{lemma}[Compact fibers over a point stratum]
\label{lem:compact-fiber-obstruction}
Consider a chain $\mathcal{M}\xrightarrow{\psi}\mathcal{P}\xrightarrow{\pi}\mathcal{B}$
with locally finite compatible refinements as in
Definition \ref{def:stratified-system}, and let
$p:\mathcal{M}\to\mathcal{M}/{\sim}$ be the quotient by
\eqref{eq:separate-relation} with the quotient topology. Suppose that
$\{q_0\}$ is a refined middle stratum, that $F=\psi^{-1}(q_0)$ is
compact and connected and meets at least two refined source strata,
and that every refined source stratum $S\subseteq F$ has a smooth
characteristic reduction inducing a strong chain. Then
$\mathcal{M}/{\sim}$ is not Hausdorff.
\end{lemma}
\begin{proof}
Compatibility makes $F$ a union of refined source strata. For such a
stratum $S$, the induced chain
$\underline S\to\{q_0\}\to\mathcal{B}_\beta$ is strong, so
Definition \ref{def:strong} gives
\begin{align}\label{eq:su3-strong-dimension}
\dim\underline S=\dim\{q_0\}+\mathrm{rank}(\pi|_{\{q_0\}})=0,
\end{align}
and $\underline S=p(S)$ is connected, hence a single point. By compactness and local finiteness, $F$ meets finitely many strata. Thus $p(F)$ is finite and connected, with at least two points because different strata are not identified. A finite Hausdorff space is discrete, a contradiction.
\end{proof}

\begin{proposition}\label{prop:unrestricted-zero}
Let $G=\mathrm{SU}(n)$, $n\ge2$, and $\xi=0$. Suppose the unrestricted chain, with its original quotient structures, admits locally finite compatible smooth stratifications refining the original strata of $\cR_0$, such that each refined stratum of $\cR_0$ has a smooth characteristic reduction inducing a strong superintegrable chain. Then $\underline{\cR}_0:=\cR_0/{\sim}$, defined by \eqref{eq:separate-relation} with the quotient topology, is not Hausdorff.
\end{proposition}
\begin{proof}
The moment equation gives $y=-x$ and $\gamma\in G_x$, so
\begin{align}\label{eq:zero-sun-chain}
 \cR_0=\{(x,\gamma):x\gamma=\gamma x\}/G,\qquad
 \cP_0\cong\mathfrak{g}/G,\qquad \pi=\mathrm{id}.
\end{align}
The fiber
\begin{align}\label{eq:zero-sun-fiber}
 F:=\psi^{-1}(0)=\{[(0,\gamma)]:\gamma\in G\}\cong G/\mathrm{Ad}G
\end{align}
is compact and connected, since $G$ is compact and connected. It contains the $n$ distinct points $m_\zeta=[(0,\zeta I_n)]$, where $\zeta^n=1$. These are exactly the points with stabilizer $G$: commuting with every element of $G$ forces $x=0$ and $\gamma=\zeta I_n$. Hence each $\{m_\zeta\}$ is an original point stratum and remains a point stratum under refinement.

Let $Q$ be the refined stratum of $\cP_0$ containing $0$. The characteristic reduction of $\{m_1\}$ is itself. Its induced map to $Q$ must be surjective by Definition \ref{def:strong}, so $Q=\{0\}$. Lemma \ref{lem:compact-fiber-obstruction}, applied to the compact connected fiber $F$ and the $n\ge2$ point strata $\{m_\zeta\}$ contained in it, proves that $\underline{\cR}_0$ is not Hausdorff. Consequently, no such locally finite stratifications make the unrestricted zero-orbit chain stratified superintegrable in the sense of Definition \ref{def:stratified-system}.
\end{proof}

\subsubsection{The restricted Liouville system and scope of the obstruction}\label{subsec:unrestricted-scope}
The simple-spectrum restriction removes this obstruction and gives an explicit Liouville chain.

\begin{proposition}
\label{prop:scope-without-spectrum}
For $\xi=0$, we have $\Omega_0=\mathcal{C}^\circ$. With the original smooth structures, the restricted chain is
\begin{align}\label{eq:zero-restricted-product}
 \mathcal{M}_0\cong\mathcal{C}^\circ\times\mathbf{H}  \xrightarrow{\mathrm{pr}_1}\mathcal{P}_0^\Omega\cong\mathcal{C}^\circ  \xrightarrow{\mathrm{id}}\Omega_0=\mathcal{C}^\circ.
\end{align}
It is a strong Liouville chain of rank $n-1$, and hence a
stratified superintegrable chain with one stratum in each space.
\end{proposition}
\begin{proof}

Every $\lambda\in\mathcal{C}$ is feasible, by taking $x=\ii\mathrm{diag}\lambda$, $\gamma=I_n$. If $\mu=0$ and the first Horn set in \eqref{hyp:spectrum} is nonempty, its matrices satisfy $A+B=0$ because their Hermitian sum has zero spectrum. Thus $\lambda[I]=\lambda[J]$ as tuples. For distinct entries this forces $I=J$, proving $\Omega_0=\mathcal{C}^\circ$.

The moment equation is $x\gamma=\gamma x$. Ordered diagonalization of simple $x$ then makes $\gamma$ diagonal, and the residual diagonal conjugation is trivial. Smooth local eigenframes therefore give \eqref{eq:zero-restricted-product} with the original smooth structures. Its form is
\begin{equation}\label{eq:zeroprod}
\omega_0=\sum_i\dd\lambda_i\wedge\dd\theta_i,\qquad
\sum_i\lambda_i=0,\quad\sum_i\dd\theta_i=0,
\end{equation}
where $\dd\theta_i=-\ii t_i^{-1}\dd t_i$ for $t=\mathrm{diag}(t_i)\in\mathbf{H}$. The trace-zero pairing is nondegenerate, the target brackets vanish, and the dimensions are $2k,k,k$. Together with the connected $k$-torus fibers, these facts give a strong Liouville chain with $p=\mathrm{id}$. The obstructing fiber \eqref{eq:zero-sun-fiber} lies over $0\notin\mathcal{C}^\circ$.

\end{proof}

Proposition \ref{prop:unrestricted-zero} concerns the unrestricted zero-orbit chain. The next subsection gives a nonzero-orbit obstruction at a repeated spectrum. Neither proves failure at every nonzero orbit or every excluded simple spectrum. At an excluded simple spectrum, Proposition \ref{prop:stabilizer} and Lemma \ref{lem:hamiltonian-torus} rule out full spectral rank on the whole fiber, but not lower-rank chains. Extensions to nonfree points, repeated spectra or wall normals without a common lattice basis require separate proofs of constant rank, Hausdorffness and compatibility of the smooth and Poisson structures.

\subsection{A compact fiber with a non-Hausdorff characteristic quotient}
\label{app:disk}\label{subsec:disk}

Fix $\xi=\ii\mathrm{diag}(a,0,-a)$, $a>0$, and put
\begin{equation}\label{eq:disk-A}
x_0=\ii\mathrm{diag}(a/3,-2a/3,a/3),\quad
y_0=-\ii\mathrm{diag}(-2a/3,a/3,a/3),\quad q_0=[(x_0,y_0)].
\end{equation}
Then $q_0\in\mathcal{P}_\xi$ has spectrum $\lambda_0=(a/3,a/3,-2a/3)$, or $(u,v)=(0,a)$. We describe the whole fiber $F=\psi^{-1}(q_0)$.

Choose $\gamma_{\rm ref}=\mathrm{diag}\bigl(\left(\begin{smallmatrix}0&1\\-1&0\end{smallmatrix}\right),1\bigr)$. Every lift of $(x_0,y_0)$ is $\gamma=h\gamma_{\rm ref}$ with $h\in G_{x_0}$, represented by a block $V\in\mathrm U(2)$ on $\mathrm{span}(e_1,e_3)$ and by $(\det V)^{-1}$ on $e_2$. The remaining conjugations form $G_{x_0,y_0}=\mathbf{H}$, acting by
\begin{equation}\label{eq:disk-gauge}
V\longmapsto\begin{pmatrix}(t_1/t_2)v_{11}&(t_1/t_3)v_{12}\\(t_3/t_2)v_{21}&v_{22}\end{pmatrix},\qquad t_1t_2t_3=1.
\end{equation}
Writing $\zeta=v_{22}$, $d=\det V$, $b=v_{12}$ gives
$V=\left(\begin{smallmatrix}d\bar\zeta&b\\-d\bar b&\zeta\end{smallmatrix}\right)$,
where $|d|=1$ and $|b|^2=1-|\zeta|^2$. The map $t\mapsto(t_1/t_2,t_1/t_3)$ is onto $(S^1)^2$: for $(u,v)$ choose $t_1^3=uv$, $t_2=t_1/u$, $t_3=t_1/v$. Thus one may set $d=1$, $b=\sqrt{1-|\zeta|^2}$, also when $b=0$, obtaining
\begin{equation}\label{eq:disk-U}
\gamma_\zeta=\begin{pmatrix}0&\bar\zeta&\sqrt{1-|\zeta|^2}\\-1&0&0\\0&-\sqrt{1-|\zeta|^2}&\zeta\end{pmatrix}\in\mathrm{SU}(3),\qquad |\zeta|\le1.
\end{equation}
Directly,
\begin{equation}\label{eq:disk-moment}
-\gamma_\zeta^{-1}x_0\gamma_\zeta=y_0,\qquad
J(x_0,\gamma_\zeta)=\ii\mathrm{diag}(a,-a,0)=:z_0.
\end{equation}
The invariant $\mathrm{tr}\gamma_\zeta=\zeta$ separates these classes. Hence the continuous parametrization $\zeta\mapsto m_\zeta=[(x_0,\gamma_\zeta)]$ is a bijection from the compact closed unit disk to $F$, and is a homeomorphism because $\mathcal R_\xi$ is Hausdorff.

Every common stabilizer lies in $G_{x_0,y_0}=\mathbf{H}$. For $|\zeta|<1$, the nonzero entries $13,21,32$ of $\gamma_\zeta$ force its diagonal entries to coincide, giving $Z_3$. At $|\zeta|=1$ the stabilizer is $H(2,1)$. The interior and boundary circle are smoothly embedded by \eqref{eq:disk-U}, with inverse $\mathrm{tr}\gamma$.

On the interior, write $\gamma_\zeta=h_\zeta\gamma_{\rm ref}$. The matrix $h_\zeta$ has an $\mathrm{SU}(2)$ block on $\mathrm{span}(e_1,e_3)$ and is the identity on $e_2$. Since $x_0$ is scalar on this block,
\[
\Theta|_F=-\mathrm{tr}(x_0\,\dd h_\zeta h_\zeta^{-1})
=-\ii(a/3)\mathrm{tr}(\dd h_\zeta h_\zeta^{-1})=0.
\]
Equation \eqref{eq:disk-moment} also makes the orbit-form pullback zero. Thus \eqref{eq:reduction-form} vanishes on the interior. It vanishes on the boundary for dimensional reasons.

The two original orbit-type intersections are
\begin{equation}\label{eq:disk-primary-intersections}
S_{\mathrm{int}}=\{m_\zeta:|\zeta|<1\}
=F\cap(\mathcal R_\xi)_{(Z_3)},\qquad
S_\partial=\{m_\zeta:|\zeta|=1\}
=F\cap(\mathcal R_\xi)_{(H(2,1))}.
\end{equation}
They are connected manifolds of dimensions two and one, respectively,
with zero restricted forms. If they are refined source strata, each
separately reduces to a point and gives the strong rank-zero chain
\begin{equation}\label{eq:disk-reduced-chains}
\{\underline m_i\}\xrightarrow{\underline\psi}\{q_0\}
\xrightarrow{\pi}\{\lambda_0\},\qquad i\in\{\mathrm{int},\partial\}.
\end{equation}
For this two-piece refinement, the interior is one characteristic leaf and is not closed:
\begin{equation}\label{eq:disk-incidence}
\overline{S_{\mathrm{int}}}=F,\qquad
m_{1-1/j}\longrightarrow m_1\in S_\partial.
\end{equation}
Its quotient class is not closed, so this characteristic quotient is
not $T_1$.

\begin{proposition}\label{prop:quotient-obstruction}
Let the fixed-orbit chain for
$\xi=\ii\mathrm{diag}(a,0,-a)$, $a>0$, admit locally finite compatible
smooth refinements which refine the original source orbit types.
Assume that $\{q_0\}$, with $q_0$ as in \eqref{eq:disk-A}, is a
refined middle stratum, and that every refined source stratum
contained in $F=\psi^{-1}(q_0)$ has a smooth characteristic reduction
inducing a strong chain. Then the characteristic quotient
\eqref{eq:separate-relation}, with its quotient topology, is not
Hausdorff.
\end{proposition}
\begin{proof}
The fiber $F$ is compact and connected and meets both original
orbit types in \eqref{eq:disk-primary-intersections}. It therefore
meets at least two refined source strata. Lemma
\ref{lem:compact-fiber-obstruction} applies.
\end{proof}
The assumption that $\{q_0\}$ is a refined middle stratum is part of
this statement. Unlike the zero-orbit argument, the disk calculation
does not force that assumption for every possible refinement.

\section*{Acknowledgements}
The authors thank Zhuo Chen, Kai Jiang, Nicolai Reshetikhin, and Husileng Xiao for discussions at an early stage of this work. The authors also thank OpenAI’s ChatGPT‑6 for helpful suggestions regarding Condition (H). The authors are responsible for the correctness of the results.
 
\bibliographystyle{unsrt}
\bibliography{bibliography}

\end{document}